\documentclass[preprint,12pt]{elsarticle}

\usepackage{amssymb}
\usepackage{amsmath}
\usepackage{amsthm}

\newtheorem{theorem}{Theorem}
\newtheorem{lemma}{Lemma}
\newtheorem{proposition}{Proposition}
\newtheorem{corollary}{Corollary}

\newtheorem{remark}{Remark}
\newtheorem{assumption}{Assumption}

\journal{Computer Methods in Applied Mechanics and Engineering}

\usepackage{amssymb}
\usepackage{mathrsfs}

\newcommand{\eq}{:=}
\newcommand{\osw}{I_{\mathrm{av}}^h}

\newcommand{\tr} {\mathrm{tr}}
\newcommand{\sym} {\mathrm{sym}}
\newcommand{\skw} {\mathrm{skw}}

\newcommand{\bsigma}{\boldsymbol{\sigma}}
\newcommand{\btau}{\boldsymbol{\tau}}
\newcommand{\brho}{\boldsymbol{\rho}}
\newcommand{\bff}{\boldsymbol{f}}
\newcommand{\bzero}{\boldsymbol{0}}
\newcommand{\bzeta}{\boldsymbol{\zeta}}

\DeclareMathOperator*{\argmin}{arg\,min}

\newcommand{\ddiv}{\operatorname{div}}

\newcommand{\jmp}[1]{[\![#1]\!]}

\newcommand{\BE}{\boldsymbol E}

\newcommand{\bg}{\boldsymbol g}

\newcommand{\bn}{\boldsymbol n}

\newcommand{\br}{\boldsymbol r}

\newcommand{\bu}{\boldsymbol u}
\newcommand{\bv}{\boldsymbol v}
\newcommand{\bw}{\boldsymbol w}
\newcommand{\bx}{\boldsymbol x}

\newcommand{\bz}{\boldsymbol z}

\newcommand{\CA}{\mathcal A}
\newcommand{\CB}{\mathcal B}
\newcommand{\CC}{\mathcal C}

\newcommand{\CE}{\mathcal E}

\newcommand{\CM}{\mathcal M}

\newcommand{\CO}{\mathcal O}

\newcommand{\CT}{\mathcal T}

\newcommand{\T}{\mathsf{T}}

\usepackage{xcolor}

\usepackage{multirow}

\usepackage{tikz}
\usetikzlibrary{arrows.meta, shapes.geometric, positioning, patterns}

\usepackage{placeins}

\begin{document}

\begin{frontmatter}



\title{A posteriori error estimates based on local residual minimization for weakly symmetric mixed finite element discretizations of linear elasticity}


\author{Ernesto Cáceres} \ead{ernesto.caceres@usach.cl}
\author{Patrick Vega\corref{cor1}}
\ead{patrick.vega@usach.cl}
\cortext[cor1]{Corresponding author}

\affiliation{
    organization={Department of Mathematics and Computer Science,
                  Universidad de Santiago de Chile},
    addressline={Las Sophoras 165, Estación Central},
    city={Santiago},
    postcode={9170020},
    country={Chile}
}

\affiliation{
    organization={Computational Heat and Fluid Flow Laboratory,
                  Universidad de Santiago de Chile},
    addressline={Av. Víctor Jara 3659, Estación Central},
    city={Santiago},
    postcode={9170124},
    country={Chile}
}

%% Abstract
\begin{abstract}
We introduce an a posteriori error estimator for weakly symmetric mixed finite element discretizations of linear elasticity, based on local residual minimization in discrete dual norms. Relying on existing a priori error estimates for four families of weakly symmetric mixed methods, we start from a local postprocessing scheme for the displacement, enlarge its test space, and minimize the resulting residual in the corresponding discrete dual norm. Here the scheme is driven by the full displacement gradient, rather than its symmetric part alone, so as to incorporate the discrete rotation multiplier. Being entirely local, this construction yields, at negligible additional cost, a residual field that serves as the main ingredient of an a posteriori error indicator. This indicator carries an additional term, with no counterpart in the exactly symmetric case, measuring the extent to which the discrete stress fails to satisfy angular momentum balance pointwise. We establish a local efficiency estimate, with an explicit constant independent of $\lambda$, and a reliability estimate in the compressible regime, with the dependence on the Lam\'e parameters tracked explicitly. By adapting a robust stability property of the continuous problem, we also obtain an alternative reliability estimate, with a constant that remains bounded as $\lambda\to\infty$. A single indicator and a single comparison norm thus serve both regimes. Numerical experiments confirm the theoretical findings.
\end{abstract}


% \begin{highlights}
% \item A posteriori estimator for weakly symmetric elasticity via residual minimization.
% \item Residual minimization on an enlarged test space of superconvergent postprocessing.
% \item Reliability and efficiency shown in both compressible and incompressible regimes.
% \item The framework covers four existing families of weakly symmetric mixed elements.
% \end{highlights}

%% Keywords

\end{frontmatter}



\section{Introduction}

Numerical simulation of elastic structures requires accurate approximation not only of the displacement field, but also of the stress tensor itself, since stress is the quantity of direct engineering interest for design and failure assessment. Mixed finite element methods, see~\cite{BoffiBrezziFortin2013} for a general reference, provide this direct stress approximation while remaining robust as the material approaches incompressibility, a regime in which standard displacement-based discretizations suffer from locking~\cite{BabuskaSuri1992} and become unreliable. Combined with adaptive mesh refinement driven by a reliable a posteriori error indicator, they are an attractive choice for practical elasticity computations.

A key physical requirement on the stress tensor is symmetry, a consequence of the balance of angular momentum. Constructing discrete stress spaces that satisfy this symmetry pointwise, as pursued by Watwood and Hartz~\cite{WatwoodHartz1968}, Johnson and Mercier~\cite{JohnsonMercier1978}, extended to higher dimensions in~\cite{GopalakrishnanGuzmanLee2026}, Arnold and Winther~\cite{ArnoldWinther2002}, and several later families~\cite{ArnoldDouglasGupta1984,GuzmanNeilan2014,ArnoldAwanouWinther2008}, is a well-known implementation challenge; see~\cite{CarstensenEigelGedicke2011} for a computational comparison. An alternative, dating back to Fraeijs de Veubeke~\cite{FraeijsDeVeubeke1977}, is to relax this requirement, enforcing symmetry only weakly through an additional Lagrange multiplier representing the infinitesimal rotation. This trade-off, an extra unknown in exchange for simpler, more readily implementable discrete stress spaces, has led to several practical element families~\cite{AmaraThomas1979,ArnoldBrezziDouglas1984,Stenberg1988,Cockburn2010,Gopalakrishnan2012,Arnold2007}, and is the setting adopted in this paper.

A priori and a posteriori error analyses for both the exactly and weakly symmetric approaches are given, respectively, in~\cite{LS2023} and~\cite{LS2024}, the latter being our starting point. A practical feature of mixed methods, going back to Arnold and Brezzi~\cite{ArnoldBrezzi1985}, is that the displacement approximation can be improved, once the discrete problem has been solved, by an elementwise postprocessing step at negligible additional cost. For second-order elliptic problems, a scheme of this type was proposed by Stenberg~\cite{Stenberg1991}; for mixed methods in elasticity, the idea dates back to Stenberg~\cite{Stenberg1988}.

This paper makes three contributions, aimed at turning this postprocessing step into an error indicator for adaptive mesh refinement. First, we enlarge the local test space of the postprocessing scheme and minimize the resulting residual in the corresponding discrete dual norm, following the methodology introduced by Muga, Rojas, and Vega~\cite{MugaRojasVega2023} for the scheme of~\cite{Stenberg1991} in scalar diffusion problems. Here, following~\cite{LS2024}, the scheme is formulated in terms of the full displacement gradient, rather than its symmetric part alone, so as to incorporate the discrete rotation on its right-hand side. Second, this construction gives, as a by-product and at negligible additional cost, a residual field that we use to build an a posteriori error indicator, with a reliability estimate in the standard compressible regime and a local efficiency estimate that holds uniformly in $\lambda$. Alongside this residual field, the indicator measures how well the gradient of the postprocessed displacement matches the discrete stress and rotation, how far this displacement is from being conforming, and, through an asymmetry term, how far the discrete stress is from satisfying angular momentum balance. Its analysis requires separating the symmetric and skew-symmetric parts of the error, so that the stress, the rotation, and the postprocessed displacement are controlled simultaneously. Third, adapting the argument of~\cite[Theorem~8]{LS2024}, we show that the same indicator also satisfies a reliability estimate with a constant that remains bounded as $\lambda\to\infty$.

\subsection{Related work}

The general theory of a posteriori error estimation is well established; see~\cite{AinsworthOden2000,Verfurth2013}. For mixed methods for second-order elliptic problems, Lovadina and Stenberg~\cite{LovadinaStenberg2006} showed that the superconvergent postprocessing of Stenberg~\cite{Stenberg1991} yields a reliable and efficient a posteriori error estimator; in elasticity, where the postprocessing dates back to Stenberg~\cite{Stenberg1988}, postprocessed displacements have since entered a posteriori estimators for both weakly and exactly symmetric methods, alongside other approaches. For weakly symmetric stresses, Carstensen and Dolzmann~\cite{CarstensenDolzmann1998} establish a residual-based estimator in two dimensions, with constants independent of $\lambda$, and Lonsing and Verf\"urth~\cite{LonsingVerfurth2004} develop a residual estimator and estimators based on auxiliary local problems for the PEERS~\cite{ArnoldBrezziDouglas1984} and BDMS~\cite{Stenberg1988} elements. Closer to the present approach, Kim~\cite{Kim2011} postprocesses the displacement through local problems driven by the full gradient, with datum $\CC\bsigma_h+\brho_h$, and uses an averaged version of the result to obtain, for the Arnold--Falk--Winther (AFW)~\cite{Arnold2007} and BDMS elements in two dimensions, a guaranteed upper bound for the stress error in the energy norm induced by $\CC$, robust in $\lambda$ by means of computable Korn and inf-sup constants, together with local lower bounds on element patches. Lederer and Stenberg~\cite{LS2024} derive estimates based on the Prager--Synge hypercircle principle, together with an alternative estimate with a constant that remains bounded as $\lambda\to\infty$. For exactly symmetric stresses, Carstensen and Gedicke~\cite{CarstensenGedicke2016} develop residual-based estimators for the Arnold--Winther element that require an efficient approximation of the skew-symmetric part of the gradient, obtained from a postprocessed displacement; Carstensen, Gallistl, and Gedicke~\cite{CarstensenGallistlGedicke2019} later remove this requirement with an explicit estimator depending only on the symmetric stress. Lederer and Stenberg~\cite{LS2023} give the exactly symmetric counterpart of the analysis in~\cite{LS2024}.

The residual-minimization methodology used here is closely related to discontinuous Petrov--Galerkin (DPG) methods~\cite{DemkowiczGopalakrishnan2011,CarstensenDemkowiczGopalakrishnan2014,Demkowicz_Gopalakrishnan_2025}, which have also been applied to several variational formulations of linear elasticity by Keith, Fuentes, and Demkowicz~\cite{KeithFuentesDemkowicz2016} and Fuentes, Keith, Demkowicz, and Le Tallec~\cite{FuentesKeithDemkowiczLeTallec2017}. Closer to the present work, this residual-minimization framework was combined with the postprocessing technique of Stenberg~\cite{Stenberg1991} by Muga, Rojas, and Vega~\cite{MugaRojasVega2023} for scalar diffusion problems, and later extended by Camargo, Rojas, and Vega~\cite{CamargoRojasVega2025} to hybridizable discontinuous Galerkin discretizations of the Helmholtz equation, and by the present authors, in a companion paper~\cite{CaceresVega}, to exactly symmetric linear elasticity. The present work extends this methodology to linear elasticity with weakly symmetric stresses, requiring a treatment of the coupled symmetric and skew-symmetric parts of the gradient with no counterpart in these earlier applications.

\subsection{Outline of the paper}

The remainder of this paper is organized as follows.  Section~\ref{sec:governing_eqns} presents the model problem, the discretization by weakly symmetric mixed finite elements, the associated a priori error estimate, and the postprocessing scheme of Stenberg, adapted to the present setting. Section~\ref{sec:residual-min} formulates a local residual minimization problem built on this postprocessing scheme, with an enlarged test space, and relates its minimizer to the postprocessed displacement. Section~\ref{sec:aposteriori} introduces an a posteriori error indicator based on the resulting residual and establishes its local efficiency together with two reliability estimates, the second of them robust as $\lambda\to\infty$. Numerical experiments illustrating these theoretical findings are reported in Section~\ref{sec:numerical_results}.

\section{Governing equations with weakly enforced stress symmetry}\label{sec:governing_eqns}

\subsection{The elasticity problem}

Let $\Omega\subset\mathbb{R}^d$ ($d\in\{2,3\}$) be a bounded polygonal or polyhedral domain, with boundary $\Gamma\eq\partial\Omega$ split into two disjoint parts $\Gamma_D$ and $\Gamma_N$, with $|\Gamma_D|>0$. The physical unknowns are the displacement vector $\bu=(u_1,\ldots,u_d):\Omega\to\mathbb{R}^d$ and the stress tensor $\bsigma=\{\sigma_{ij}\}_{i,j=1}^d:\Omega\to\mathbb{R}^{d\times d}$, which is \emph{not} assumed symmetric, in contrast to the exactly symmetric case. The loading consists of a body force $\bff$ and a traction $\bg$ prescribed on the boundary part $\Gamma_N$; on the complementary part $\Gamma_D$, homogeneous Dirichlet conditions are imposed on the displacement.

The balances of linear and angular momentum give rise, respectively, to the equations
\begin{subequations}
\begin{align}
\operatorname{div}\bsigma + \bff &= \bzero \quad \text{in }\Omega, 
&&\text{(linear momentum balance)} \label{eq:linmom}\\
\bsigma - \bsigma^{\T} &= \bzero \quad \text{in }\Omega. 
&&\text{(angular momentum balance)} \label{eq:angmom}
\end{align}
\end{subequations}

In the static setting considered here, \eqref{eq:linmom} and \eqref{eq:angmom} express translational and rotational equilibrium, respectively; in the absence of body couples, \eqref{eq:angmom} is the usual physical requirement that the stress tensor be symmetric.

The deformation is characterized by the infinitesimal strain tensor

\begin{equation}
\varepsilon(\bu) = \frac{1}{2}\big(\nabla \bu + (\nabla \bu)^{\T}\big).
\label{eq:strain}
\end{equation}

For an isotropic material, the constitutive relation is linear, $\CC\bsigma = \varepsilon(\bu)$, with compliance tensor $\CC\btau$ given, depending on the regime under consideration, in Table~\ref{tab:compliance}.

\begin{table}[h]
\centering
\begin{tabular}{ll}
\hline
\textbf{\small Regime} & \textbf{\small Compliance tensor} \\
\hline
\small Plane strain ($d=2$) & \small
\multirow{2}{*}{$\CC\btau = \dfrac{1}{2\mu}\Big(\btau - \dfrac{\lambda}{2\mu+d\lambda}\tr(\btau)I\Big)$} \\
3D ($d=3$) & \\[8pt]
\small Plane stress ($d=2$) & \small
$\CC\btau = \dfrac{1}{2\mu}\Big(\btau - \dfrac{\lambda^*}{2\mu+2\lambda^*}\tr(\btau)I\Big)$, \quad \small$\lambda^* = \dfrac{2\mu\lambda}{\lambda+2\mu}$ \\
\hline
\end{tabular}
\caption{Compliance tensor $\CC\btau$ for the regimes considered. Here $\mu,\lambda$ denote the Lam\'e parameters; for plane stress, all results below hold with $\lambda$ replaced by $\lambda^*$.}
\label{tab:compliance}
\end{table}

The inverse of the compliance tensor, the elasticity tensor, is denoted $\CA\eq\CC^{-1}$, given by $\CA\btau = 2\mu\btau + \lambda\,\tr(\btau)I$ (respectively with $\lambda^*$ in place of $\lambda$ for plane stress). In the limit $\lambda\to\infty$ the material becomes incompressible, characterized by the constraint $\operatorname{div}\bu = 0$; as $\lambda^*<2\mu$, this limit concerns only plane strain and the three-dimensional case.

Since the compliance tensor is itself symmetric, symmetry of the stress is implicitly built into this relation: the traditional form of the equations of elasticity reads
\begin{subequations}
\begin{align}
\CC\bsigma - \varepsilon(\bu) &= \bzero \quad \text{in }\Omega, 
&&\text{(constitutive equation)} \label{eq:constitutive}\\
\operatorname{div}\bsigma + \bff &= \bzero \quad \text{in }\Omega, 
&&\text{(linear momentum balance)} \label{eq:equilibrium}\\
\bu &= \bzero \quad \text{on }\Gamma_D, 
&&\text{(clamping condition)} \label{eq:dirichlet}\\
\bsigma \bn &= \bg \quad \text{on }\Gamma_N. 
&&\text{(traction condition)} \label{eq:traction}
\end{align}
\end{subequations}

%\subsection{Why an additional unknown is needed}

Using this formulation directly as the basis for a mixed finite element method, however, leads to methods that are cumbersome and costly to implement, owing precisely to the difficulty of constructing discrete spaces that satisfy stress symmetry exactly. This observation led Fraeijs de Veubeke to propose mixed methods in which rotational equilibrium \eqref{eq:angmom} is treated as an \emph{independent equation}, rather than being enforced exactly at the level of the discrete space. This relaxation requires introducing an additional unknown: the skew-symmetric infinitesimal rotation
\begin{equation}
\brho = \frac{1}{2}\big(\nabla \bu - (\nabla \bu)^{\T}\big).
\label{eq:rho-def}
\end{equation}

With this, the complete system of equations is rewritten as
\begin{subequations}\label{eq:full_system}
\begin{align}
\CC\bsigma + \brho - \nabla \bu &= \bzero \quad \text{in }\Omega, 
&&\text{(constitutive equation)} \label{eq:2.15}\\
\bsigma - \bsigma^{\T} &= \bzero \quad \text{in }\Omega, 
&&\text{(angular momentum balance)} \label{eq:2.16}\\
\operatorname{div}\bsigma + \bff &= \bzero \quad \text{in }\Omega, 
&&\text{(linear momentum balance)} \label{eq:2.17}\\
\bu &= \bzero \quad \text{on }\Gamma_D, 
&&\text{(clamping condition)} \label{eq:2.18}\\
\bsigma \bn &= \bg \quad \text{on }\Gamma_N. 
&&\text{(traction condition)} \label{eq:2.19}
\end{align}
\end{subequations}

For any tensor $\btau\in[L^2(\Omega)]^{d\times d}$, we write $\sym(\btau)\eq\frac{1}{2}(\btau+\btau^{\T})$ and $\skw(\btau)\eq\frac{1}{2}(\btau-\btau^{\T})$ for its symmetric and skew-symmetric parts, so that $\btau=\sym(\btau)+\skw(\btau)$. We denote by \[[L^2(\Omega)]^{d\times d}_{\skw}\eq\{\btau\in[L^2(\Omega)]^{d\times d}\mid \btau^{\T}=-\btau\}\] the subspace of skew-symmetric tensor fields. In what follows we use the notation
\begin{equation}
\omega(\bu) \eq \skw(\nabla\bu) = \frac{1}{2}\big(\nabla \bu - (\nabla \bu)^{\T}\big),
\label{eq:omega-def}
\end{equation}
so that $\brho=\omega(\bu)$, and $\varepsilon(\bu)=\sym(\nabla\bu)$. Since $\sym(\btau)$ and $\skw(\btau)$ are orthogonal in the $L^2$-inner product, $\|\btau\|_0^2=\|\sym(\btau)\|_0^2+\|\skw(\btau)\|_0^2$, and in particular
\begin{equation}
\|\sym(\btau)\|_0 \le \|\btau\|_0 \quad\text{and}\quad \|\skw(\btau)\|_0 \le \|\btau\|_0 \qquad \forall\btau\in[L^2(\Omega)]^{d\times d}.
\label{eq:sym-skw-bound}
\end{equation}

\begin{remark}[Symmetric/skew-symmetric splitting]\label{rem:splitting}
Since $\CC\bsigma$ is symmetric and $\brho$ is skew-symmetric, the symmetric and skew-symmetric parts of equation~\eqref{eq:2.15} decouple exactly: its symmetric part is precisely the classical constitutive relation~\eqref{eq:constitutive}, while its skew-symmetric part is precisely the identity $\brho=\omega(\bu)$ noted above, which identifies $\brho$ with the infinitesimal rotation.
\end{remark}

The Lam\'e parameters $\mu,\lambda$ appear throughout the error estimates in this work, so we fix here the notation used to track their dependence: we write $A\lesssim B$ (and $A\gtrsim B$) if there exists a positive constant $C$, independent of the mesh size $h$ (but possibly depending on $\Omega$, $\Gamma_D$, and the shape regularity of the mesh family), such that $A\le CB$ (respectively, $A\ge CB$). Unless stated otherwise, $C$ may depend on $\mu,\lambda$; whenever an estimate holds with $C$ independent of these parameters as well, this is indicated explicitly in the statement.

\subsection{Weakly symmetric mixed finite element methods}

The system~\eqref{eq:full_system} admits a dual-mixed variational formulation: find $(\bsigma,\brho,\bu)\in H_g(\ddiv;\Omega)\times[L^2(\Omega)]^{d\times d}_{\skw}\times[L^2(\Omega)]^d$ such that
\begin{align}\label{eq:dual-variational}
\CM(\bsigma,\brho,\bu;\btau,\bzeta,\bv) + (\bff,\bv) &= 0\\ &\forall(\btau,\bzeta,\bv)\in H_0(\ddiv;\Omega)\times[L^2(\Omega)]^{d\times d}_{\skw}\times[L^2(\Omega)]^d,\nonumber
\end{align}
with
\begin{align*}
    \CM(\bsigma,\brho,\bu;\btau,\bzeta,\bv)&\eq(\CC\bsigma,\btau)+b(\btau;\brho,\bu)+b(\bsigma;\bzeta,\bv),\\
    b(\btau;\bzeta,\bv)&\eq(\ddiv\btau,\bv)+(\btau,\bzeta),
\end{align*}
and
\begin{align}
H_g(\ddiv;\Omega)&\eq\{\btau\in H(\ddiv;\Omega)\mid \btau\bn|_{\Gamma_N}=\bg\},\\
H_0(\ddiv;\Omega)&\eq\{\btau\in H(\ddiv;\Omega)\mid \btau\bn|_{\Gamma_N}=\bzero\}.
\end{align}

Following~\cite{LS2024}, we consider finite element discretizations based on the variational formulation~\eqref{eq:dual-variational}, posed over subspaces $S_h\subset H(\ddiv;\Omega)$, $R_h\subset[L^2(\Omega)]^{d\times d}_{\skw}$, and $V_h\subset[L^2(\Omega)]^d$, on a triangular or tetrahedral mesh $\CT_h$, belonging to a shape-regular family of conforming triangulations. Unified a priori error analyses of weakly symmetric mixed methods, with estimates in $L^2$ norms, are available~\cite{Guzman2010,Lee2016}: both cover the Cockburn--Gopalakrishnan--Guzm\'an (CGG)~\cite{Cockburn2010}, Stenberg--Gopalakrishnan--Guzm\'an (SGG)~\cite{Gopalakrishnan2012,Stenberg1988}, and AFW families, and the framework of~\cite{Lee2016} also covers elements with continuous rotations, such as PEERS~\cite{ArnoldBrezziDouglas1984}. Here we rely instead on the analysis of~\cite{LS2024}, which provides a priori estimates in the mesh-dependent norm~\eqref{eq:broken-norm}, the one required by the postprocessing analysis below, uniformly valid in the incompressible limit.

The displacement space is indexed by a polynomial degree $k\geq 1$ and is common to every family considered:
\begin{equation}
V_h \eq \{\bv\in[L^2(\Omega)]^d \mid \bv|_K \in [P_{k-1}(K)]^d\ \ \forall K\in\CT_h\}.
\label{eq:Vh}
\end{equation}

Four families of stress/rotation pairs $(S_h,R_h)$ achieve stability, differing in how the rotation multiplier is discretized and in whether stress-space bubble functions are required to compensate for it. Table~\ref{tab:families} summarizes their defining polynomial degrees and the resulting convergence rate for the stress, the quantity that most clearly distinguishes the four families; the corresponding rates for the rotation and the (superconvergent) displacement, which follow the same pattern across all families, are discussed in Remark~\ref{rem:superconvergence} below. Notably, CGG and SGG attain the optimal rate for the stress, whereas AFW does not: its polynomial choices for the rotation and stress spaces are not in balance, and the resulting rate, $\mathcal{O}(h^k)$, falls short of the interpolation estimate $\inf_{\btau\in S_h}\|\bsigma-\btau\|_0=\mathcal{O}(h^{k+1})$ available for its stress space. The RAFW family is a modification of AFW, reducing the polynomial order of the normal trace of the stress on element interfaces, which lowers the number of degrees of freedom while retaining the same (suboptimal) accuracy at reduced computational cost.

\begin{table}[h]
\centering
\small
\renewcommand{\arraystretch}{1.5}
\setlength{\tabcolsep}{8pt}
\begin{tabular}{lccl}
\hline
Family & $R_h$ (degree) & $S_h$ (degree) & Convergence rate \\
\hline
CGG~\cite{Cockburn2010} & $k-1$ & $k-1$ + bubble & $\|\bsigma-\bsigma_h\|_0 \lesssim h^k$ \\
SGG~\cite{Gopalakrishnan2012,Stenberg1988} & $k$ & $k$ + bubble & $\|\bsigma-\bsigma_h\|_0 \lesssim h^{k+1}$ \\
AFW~\cite{Arnold2007} & $k-1$ & $k$ (no bubble) & $\|\bsigma-\bsigma_h\|_0 \lesssim h^k$ \\
RAFW~\cite{LS2024} & $k-1$ & $k$, reduced normal trace & $\|\bsigma-\bsigma_h\|_0 \lesssim h^k$ \\
\hline
\end{tabular}
\caption{Weakly symmetric families considered, with rotation space $R_h$ and stress space $S_h$ indexed by their local polynomial degree. The CGG and RAFW families are defined for $k\ge2$, and SGG and AFW for $k\ge1$~\cite{LS2024}.}
\label{tab:families}
\end{table}

The stability of the resulting bubble spaces relies on a divergence-free construction obtained from a curl-curl operator applied to a scalar (or tensor, for $d=3$) bubble function, chosen so that the bubbles vanish, together with their normal trace, on the element boundary; we omit the explicit construction here, referring the reader to~\cite{LS2024} for full detail, since our focus lies elsewhere.

We denote by $\CE_h$ the set of edges (or faces) of $\CT_h$, assumed compatible with the partition $\Gamma=\Gamma_D\cup\Gamma_N$, and by $\CE_h^0$, $\CE_h^D$, and $\CE_h^N$ its subsets of interior edges (or faces) and of those lying on $\Gamma_D$ and on $\Gamma_N$, respectively; for $K\in\CT_h$, $\CE_K^0$ and $\CE_K^D$ denote the edges (or faces) of $K$ belonging to $\CE_h^0$ and $\CE_h^D$. We write $h_K$ and $h_E$ for the diameters of $K\in\CT_h$ and $E\in\CE_h$, and, for $E\in\CE_h^0$ shared by $K^+,K^-\in\CT_h$, $\jmp{\bv}|_E\eq\bv|_{K^+}-\bv|_{K^-}$ for the jump of $\bv$ across $E$. For each $E\in\CE_h^N$, let $Q_E:[L^2(E)]^d\to[P_m(E)]^d$ denote the $L^2$-projection ($m=k$ for SGG and AFW, $m=k-1$ otherwise), and let $Q_h$ be defined edgewise/facewise by $Q_h|_E\eq Q_E$. The discrete trial and test spaces are then
\begin{subequations}
\begin{align}
S_h^g &\eq \{\btau\in S_h \mid \btau\bn = Q_h\bg\ \text{on }\Gamma_N\}, \label{eq:Shg}\\
S_h^0 &\eq \{\btau\in S_h \mid \btau\bn = \bzero\ \text{on }\Gamma_N\}. \label{eq:Sh0}
\end{align}
\end{subequations}
The discrete problem reads: find $(\bsigma_h,\brho_h,\bu_h)\in S_h^g\times R_h\times V_h$ such that
\begin{equation}\label{eq:3.3}
\CM(\bsigma_h,\brho_h,\bu_h;\btau,\bzeta,\bv) + (\bff,\bv) = 0 \qquad \forall(\btau,\bzeta,\bv)\in S_h^0\times R_h\times V_h,
\end{equation}
% with 
% \begin{align*}
%     \CM(\bsigma,\brho,\bu;\btau,\bzeta,\bv)&\eq(\CC\bsigma,\btau)+b(\btau;\brho,\bu)+b(\bsigma;\bzeta,\bv),\\ b(\btau;\bzeta,\bv)&\eq(\ddiv\btau,\bv)+(\btau,\bzeta).
% \end{align*}
with $\CM$ and $b$ as in~\eqref{eq:dual-variational}.

We also denote by $H^1(\CT_h)\eq\{v\in L^2(\Omega)\mid v|_K\in H^1(K)\quad\forall K\in\CT_h\}$ the broken Sobolev space. For $\bv\in[H^1(\CT_h)]^d$, the operators $\nabla$, $\varepsilon$, and $\omega$ are understood elementwise, and so are the corresponding $L^2(\Omega)$-norms. For the error analysis we use the broken energy norm for the displacement,
\begin{equation}
\|\bv\|_h^2 \eq \sum_{K\in\CT_h}\|\varepsilon(\bv)\|_K^2 + \sum_{E\in\CE_h^0} h_E^{-1}\|\jmp{\bv}\|_E^2 + \sum_{E\in\CE_h^D} h_E^{-1}\|\bv\|_E^2, \qquad \bv\in[H^1(\CT_h)]^d.
\label{eq:broken-norm}
\end{equation}

Stability of~\eqref{eq:3.3} with respect to this norm is established, family by family, in~\cite[Section~3]{LS2024}; as our focus lies elsewhere, we omit this analysis entirely and recall only the resulting a priori error estimate.

\begin{theorem}[A priori error estimate]
\label{thm:apriori}
Let $(\bsigma,\brho,\bu)$ solve~\eqref{eq:dual-variational}, and let\\ $(\bsigma_h,\brho_h,\bu_h)\in S_h^g\times R_h\times V_h$ solve~\eqref{eq:3.3}. For a constant independent of $\mu$ and $\lambda$, it holds
\begin{align}
&\mu^{-1/2}\|\bsigma-\bsigma_h\|_0 + \mu^{1/2}\big(\|P_h\bu-\bu_h\|_h + \|\brho-\brho_h\|_0\big)\label{eq:apriori}\\
&\lesssim \mu^{-1/2}\Bigg(\inf_{\btau\in S_h^g}\|\bsigma-\btau\|_0 + \Bigg(\sum_K h_K^2\|\bff-\bff_h\|_K^2\Bigg)^{\!1/2\,}\Bigg) + \inf_{\bzeta\in R_h}\mu^{1/2}\|\brho-\bzeta\|_0,
\nonumber
\end{align}
where $P_h:[L^2(\Omega)]^d\to V_h$ is the $L^2$-projection and $\bff_h$ is any piecewise polynomial approximation of $\bff$.
\end{theorem}
\begin{proof}
See~\cite[Theorem~4]{LS2024}.
\end{proof}

\noindent Compared with the exactly symmetric case, the right-hand side above carries an additional term, $\inf_{\bzeta\in R_h}\mu^{1/2}\|\brho-\bzeta\|_0$, reflecting the approximability of the rotation multiplier.

\begin{remark}[Superconvergence of the projected displacement error]\label{rem:superconvergence}
For a sufficiently smooth exact solution, Theorem~\ref{thm:apriori} yields the rate $\|\bsigma-\bsigma_h\|_0=\CO(h^{k+1})$ for SGG and $\CO(h^k)$ for CGG, AFW, and RAFW, and the same rate for $\|\brho-\brho_h\|_0$. The estimate for $\|P_h\bu-\bu_h\|_h$ is a superconvergence result: exploiting the divergence-preserving property $\ddiv\btau\in V_h$ for $\btau\in S_h$, one obtains $\|P_h\bu-\bu_h\|_h=\CO(h^{k+1})$ for SGG and $\CO(h^k)$ otherwise, a rate higher than that available for the unprojected error $\|\bu-\bu_h\|_h$ by a standard triangle-inequality argument combined with the approximation properties of $P_h$. This superconvergence is precisely the property exploited by the postprocessing scheme introduced in the next section.
\end{remark}

\subsection{Postprocessing of the displacement}

The idea of postprocessing the displacement in mixed methods for elasticity was originally introduced by Stenberg~\cite{Stenberg1988}, and yields a superconvergent approximation at a negligible additional computational cost, since the postprocessing is entirely local. This scheme is applicable regardless of whether the stress symmetry is imposed weakly or exactly; here, following~\cite{LS2024}, it is built from the discrete solution $(\bsigma_h,\brho_h,\bu_h)$ of~\eqref{eq:3.3}.

Let $V^l(K)\eq[P_l(K)]^d$, with $l=k$ for CGG, AFW, and RAFW, and $l=k+1$ for SGG, and let $V_h^l\eq\{\bv\in[L^2(\Omega)]^d \mid \bv|_K\in V^l(K)\ \ \forall K\in\CT_h\}$ denote its associated global space. In the original scheme of~\cite{Stenberg1988}, as in the exactly symmetric case~\cite{CaceresVega}, the postprocessing is driven by the symmetric gradient $\varepsilon(\cdot)$ paired with datum $\CC\bsigma_h$. Here, following~\cite{LS2024}, it is driven by the full gradient $\nabla(\cdot)$, paired with datum $\CC\bsigma_h+\brho_h$, which approximates $\nabla\bu$ by~\eqref{eq:2.15}; a postprocessing based on the full gradient was already used by Stenberg~\cite{Stenberg1990} for the Stokes problem, where the tensor unknown likewise approximates the full gradient of a vector field, and local problems of the same type, with the same datum, are used for a posteriori error estimation in~\cite{Kim2011}. Moreover, since $(\nabla\bv,\nabla\bz)_K=\sum_{i=1}^d(\nabla v_i,\nabla z_i)_K$ and $P_h$ acts componentwise, \eqref{eq:step1a}--\eqref{eq:step1b} decouples into $d$ scalar problems, each being the postprocessing scheme of~\cite{Stenberg1991} with datum the corresponding row of $\CC\bsigma_h+\brho_h$. This choice has two advantages in the present setting: the rotation multiplier $\brho_h$ enters the postprocessing step, and~\eqref{eq:step1a}--\eqref{eq:step1b} is uniquely solvable for every $k\ge1$, since $\|\nabla(\cdot)\|_K$ is a norm on $(I-P_h)V^l_h\big|_K$, whereas for $k=1$ this space contains nonzero infinitesimal rotations, on which $\varepsilon$ vanishes. We look for $\widetilde{\bu}_h\in V_h^l$ satisfying the following equations on each $K\in\CT_h$:
\begin{subequations}
\begin{alignat}{2}
(\nabla\widetilde{\bu}_h,\nabla\bv_K)_K &= (\CC\bsigma_h+\brho_h,\nabla\bv_K)_K \qquad && \forall \bv_K\in(I-P_h)V^l_h\big|_K, \label{eq:step1a}\\
(\widetilde{\bu}_h,\bv_K)_K &= (\bu_h,\bv_K)_K \qquad && \forall \bv_K\in P_h V_h^l\big|_K. \label{eq:step1b}
\end{alignat}
\end{subequations}

A closely related scheme for vectorial problems, with the resulting a priori estimate expressed in terms of quantities already controlled by the a priori error analysis of the underlying method, is analyzed in~\cite{ArayaSolanoVega2019} for a hybridizable discontinuous Galerkin discretization of the Brinkman problem; the same type of estimate is established for the scalar case, within a general framework illustrated with the Brezzi--Douglas--Marini element, in~\cite{MugaRojasVega2023}, where it then serves as the basis of a local residual minimization problem, the approach followed in Section~\ref{sec:residual-min} below.

For the postprocessed displacement, we measure the error in the broken norm based on the full gradient,
\begin{equation}
\|\bv\|_{h,\nabla}^2 \eq \sum_{K\in\CT_h}\|\nabla\bv\|_K^2 + \sum_{E\in\CE_h^0} h_E^{-1}\|\jmp{\bv}\|_E^2 + \sum_{E\in\CE_h^D} h_E^{-1}\|\bv\|_E^2,
\label{eq:broken-norm-grad}
\end{equation}
defined for $\bv\in[H^1(\CT_h)]^d$. Clearly $\|\bv\|_h\le\|\bv\|_{h,\nabla}$; conversely, since \[\|\bv\|^2_{L^2(\Gamma_D)}\le\max_{E\in\CE_h^D}h_E\sum_{E\in\CE_h^D}h_E^{-1}\|\bv\|_E^2,\]
Korn's inequality for piecewise $H^1$ vector fields~\cite[eq.~(1.19)]{Brenner2004}, with $\Gamma_D$ as the boundary portion, gives $\|\bv\|_{h,\nabla}\lesssim\|\bv\|_h$ for all $\bv\in V_h$.

\begin{theorem}[A priori estimate for the postprocessed displacement]
\label{thm:postproc-apriori}
Let $(\bsigma,\brho,\bu)$ solve~\eqref{eq:dual-variational}; let $(\bsigma_h,\brho_h,\bu_h)\in S_h^g\times R_h\times V_h$ solve~\eqref{eq:3.3}; and let $\widetilde{\bu}_h\in V_h^l$ be defined by~\eqref{eq:step1a}--\eqref{eq:step1b}. Then, with a constant independent of $\mu$ and $\lambda$, it holds that
\begin{equation}\label{eq:apriori_postproc}
    \|\bu-\widetilde{\bu}_h\|_{h,\nabla} \lesssim \|\bu-P_h^l\bu\|_{h,\nabla} + \mu^{-1}\|\bsigma-\bsigma_h\|_0 + \|\brho-\brho_h\|_0 + \|P_h\bu-\bu_h\|_h,
\end{equation}
where $P_h^l:[L^2(\Omega)]^d\to V_h^l$ denotes the $L^2$-projection.
\end{theorem}
\begin{proof}
The same postprocessing is analyzed in the norm~\eqref{eq:broken-norm} in~\cite[Lemma~9]{LS2024}, where the proof is omitted; here we follow the argument of~\cite{ArayaSolanoVega2019} in the case of a vanishing reaction coefficient, the only additional step being the passage from the full gradient to the symmetric gradient in the last term. Let $\boldsymbol{e}\eq\widetilde{\bu}_h-P_h^l\bu$; since $P_hP_h^l=P_h$, \eqref{eq:step1b} gives $P_h\boldsymbol{e}=\bu_h-P_h\bu$. Since $\nabla\bu=\CC\bsigma+\brho$ by~\eqref{eq:2.15}, subtracting from~\eqref{eq:step1a} gives
\begin{equation*}
    (\nabla(\widetilde{\bu}_h-\bu),\nabla\bv_K)_K=-\big(\CC(\bsigma-\bsigma_h)+(\brho-\brho_h),\nabla\bv_K\big)_K \qquad \forall\bv_K\in(I-P_h)V^l_h\big|_K.
\end{equation*}
With $\bv_K=(I-P_h)\boldsymbol{e}|_K$, and using that $\|\CC\btau\|_K\le\mu^{-1}\|\btau\|_K$ (see~\eqref{eq:C_bound} below), yields
\begin{equation*}
    \|\nabla(I-P_h)\boldsymbol{e}\|_K\le\|\nabla(\bu-P_h^l\bu)\|_K+\mu^{-1}\|\bsigma-\bsigma_h\|_K+\|\brho-\brho_h\|_K+\|\nabla(P_h\bu-\bu_h)\|_K.
\end{equation*}
Since $(I-P_h)\boldsymbol{e}$ has zero mean on each element, a trace inequality and a Poincar\'e inequality bound its jump and Dirichlet contributions by its broken gradient. Finally, the result follows from the decomposition \[\bu-\widetilde{\bu}_h=(\bu-P_h^l\bu)-(I-P_h)\boldsymbol{e}-P_h\boldsymbol{e},\] the triangle inequality, and the bound $\|P_h\bu-\bu_h\|_{h,\nabla}\lesssim\|P_h\bu-\bu_h\|_h$ recalled after~\eqref{eq:broken-norm-grad}.
\end{proof}

Theorems~\ref{thm:postproc-apriori} and~\ref{thm:apriori} together give $\|\bu-\widetilde{\bu}_h\|_{h,\nabla}=\CO(h^l)$, that is, $\CO(h^{k+1})$ for SGG and $\CO(h^k)$ otherwise, for a sufficiently smooth exact solution, the postprocessed displacement $\widetilde{\bu}_h$ is superconvergent.

\section{Postprocessing scheme via local residual minimization}
\label{sec:residual-min}

We adapt the local residual minimization methodology of \cite{MugaRojasVega2023} to the present setting: we enlarge the local test space of the postprocessing scheme~\eqref{eq:step1a}--\eqref{eq:step1b} and minimize the resulting residual in the corresponding discrete dual norm. The minimizer coincides with the postprocessed displacement, as shown below, and the construction is valuable because it produces, as a by-product and at negligible additional cost, a Riesz representative of the associated residual, the quantity that will drive the a posteriori error indicator of Section~\ref{sec:aposteriori}.

We define a new postprocessed approximation $\bw_h$ through the following local residual minimization problems: find $\bw_h\in V_h^l$ such that for each $K\in\CT_h$,

\begin{subequations}\label{eq:minres}
    \begin{align}
    \bw_h|_K &= \argmin_{\bv_K\in V^l(K)} \frac{1}{2}\|\CC\bsigma_h+\brho_h-\nabla\bv_K\|_{*,K}^2 \label{eq:minres-a}\\
    &\qquad\text{subject to}\quad (\bv_K,\boldsymbol{q})_K = (\bu_h,\boldsymbol{q})_K \quad \forall\boldsymbol{q}\in P_h V_h^l\big|_K, \label{eq:minres-b}
    \end{align}
\end{subequations}
where $\|\cdot\|_{*,K}$ is the discrete dual norm
\begin{equation}
    \|\cdot\|_{*,K} \eq \sup_{\bv_K\in(I-P_h)V^{l+1}_h|_K\setminus\{\bzero\}} \frac{(\ \cdot\ ,\nabla\bv_K)_K}{\|\nabla\bv_K\|_K}.
    \label{eq:dual-norm}
\end{equation}

That is, $\|\btau\|_{*,K}$ is the norm, in the dual of $(I-P_h)V^{l+1}_h\big|_K$ endowed with $\|\nabla(\cdot)\|_K$, of the functional $\bv_K\mapsto(\btau,\nabla\bv_K)_K$; here $\|\nabla(\cdot)\|_K$ is indeed a norm on $(I-P_h)V^{l+1}_h\big|_K$, since $V_h$ contains the piecewise constants. As is standard for minimum-residual methods, problem~\eqref{eq:minres-a}--\eqref{eq:minres-b} can be written equivalently as the following local mixed linear system~\cite{CohenDahmenWelper2012}: find $(\br_h,\bw_h)\in(I-P_h)V^{l+1}_h\times V^l_h$ such that for each $K\in\CT_h$,
\begin{subequations}
    \begin{align}
    (\nabla\br_h,\nabla\bv_K)_K + (\nabla\bw_h,\nabla\bv_K)_K &= (\CC\bsigma_h+\brho_h,\nabla\bv_K)_K \label{eq:mixed-a}\\
    &\hspace{7em} \forall\bv_K\in(I-P_h)V^{l+1}_h\big|_K, \nonumber\\[4pt]
    (\nabla\bv_K,\nabla\br_h)_K &= 0 \label{eq:mixed-b}\\
    &\hspace{7em} \forall\bv_K\in(I-P_h)V^l_h\big|_K, \nonumber\\[4pt]
    (\bw_h,\bv_K)_K &= (\bu_h,\bv_K)_K \label{eq:mixed-c}\hspace{1.25em} \forall\bv_K\in P_h V_h^l\big|_K.
    \end{align}
\end{subequations}

The next result shows that this construction reproduces the postprocessing scheme~\eqref{eq:step1a}--\eqref{eq:step1b}.
\begin{proposition}
\label{prop:equivalence}
Let $\widetilde{\bu}_h\in V_h^l$ solve~\eqref{eq:step1a}--\eqref{eq:step1b}, and let $\bw_h\in V_h^l$ solve the residual minimization problem~\eqref{eq:minres-a}--\eqref{eq:minres-b}. Then $\bw_h=\widetilde{\bu}_h$.
\end{proposition}
\begin{proof}
Since $(I-P_h)V^l_h\big|_K\subset(I-P_h)V^{l+1}_h\big|_K$, we may test~\eqref{eq:mixed-a} with $\bv_K\in(I-P_h)V^l_h\big|_K$; by~\eqref{eq:mixed-b}, the first term vanishes, so that $\bw_h$ solves~\eqref{eq:step1a}--\eqref{eq:step1b}. The result follows from the uniqueness of the solution to~\eqref{eq:step1a}--\eqref{eq:step1b}.
\end{proof}

Hereafter we denote both postprocessed approximations by $\bw_h$. As a consequence of Proposition~\ref{prop:equivalence}, the a priori estimate of Theorem~\ref{thm:postproc-apriori} carries over unchanged to $\bw_h$.

\begin{corollary}
\label{cor:apriori-residual}
The postprocessed approximation $\bw_h\in V_h^l$ obtained from~\eqref{eq:minres} satisfies the a priori error estimate~\eqref{eq:apriori_postproc} of Theorem~\ref{thm:postproc-apriori}.
\end{corollary}

\section{A posteriori error analysis}
\label{sec:aposteriori}

We now turn to the a posteriori error analysis, based on the residual pair obtained from the postprocessing scheme of Section~\ref{sec:residual-min}. The local error indicator is built entirely from quantities already available once the discrete problem~\eqref{eq:3.3} and the postprocessing scheme have been solved, through purely local computations.

\subsection{A posteriori error indicator}

For $\bv\in[H^1(\CT_h)]^d$ and $K\in\CT_h$, let
\begin{equation}
\|\bv\|_{h,\nabla,K}^2 \eq \|\nabla\bv\|_K^2 + \frac{1}{2}\sum_{E\in\CE_K^0} h_E^{-1}\|\jmp{\bv}\|_E^2 + \sum_{E\in\CE_K^D} h_E^{-1}\|\bv\|_E^2,
\label{eq:hK-norm}
\end{equation}
so that $\|\bv\|_{h,\nabla}^2=\sum_{K\in\CT_h}\|\bv\|_{h,\nabla,K}^2$.

Since the local minimization problem is driven by the full gradient, the residual pair $(\br_h,\bw_h)$ carries information about the full displacement gradient, including its skew-symmetric part. Correspondingly, besides the Riesz representative $\br_h$ of the minimized residual, the local indicator is assembled from three physically distinct contributions: a stress residual, measuring how far the postprocessed data $\CC\bsigma_h+\brho_h$ is from being an exact gradient; a jump contribution, measuring the lack of continuity of $\bw_h$ across element interfaces; and an asymmetry defect, $\|\bsigma_h-\bsigma_h^{\T}\|_K$, quantifying the extent to which the discrete stress fails to satisfy angular momentum balance~\eqref{eq:angmom} pointwise, a quantity with no counterpart in the exactly symmetric setting, since there $\bsigma_h$ is symmetric by construction. Accordingly, we define the local error indicator $\eta_K$ by
\begin{align}
\eta_K^2 \eq\ &\mu\|\nabla\br_h\|_K^2 + \mu\|\CC\bsigma_h+\brho_h-\nabla\bw_h\|_K^2 + \mu^{-1}\|\bsigma_h-\bsigma_h^{\T}\|_K^2\notag\\
&+ \frac{\mu}{2}\sum_{E\in\CE_K^0}h_E^{-1}\|\jmp{\bw_h}\|_E^2 + \mu\sum_{E\in\CE_K^D}h_E^{-1}\|\bw_h\|_E^2,
\label{eq:etaK}
\end{align}
and the global indicator $\eta_h\eq\big(\sum_{K\in\CT_h}\eta_K^2\big)^{1/2}$.

To establish reliability with an explicit constant, we introduce the following saturation assumption, following the residual-minimization framework of~\cite{MugaRojasVega2023}, comparing $\bw_h$ against an auxiliary, higher-order approximation $\bz_h$.

Let $\bz_h\in V_h^{l+1}$ solve, for each $K\in\CT_h$,
\begin{subequations}
\begin{alignat}{2}
(\nabla\bz_h,\nabla\bv_K)_K &= (\CC\bsigma_h+\brho_h,\nabla\bv_K)_K \qquad && \forall\bv_K\in(I-P_h)V^{l+1}_h\big|_K, \label{eq:sat-a}\\
(\bz_h,\bv_K)_K &= (\bu_h,\bv_K)_K \qquad && \forall\bv_K\in P_h V_h^l\big|_K. \label{eq:sat-b}
\end{alignat}
\end{subequations}

\begin{assumption}[Saturation]
\label{ass:saturation}
There exists an $h$-independent constant $\delta\in[0,1)$ such that
\begin{equation}
\|\nabla(\bu-\bz_h)\|_0 \leq \delta\,\|\nabla(\bu-\bw_h)\|_0.
\label{eq:saturation}
\end{equation}
\end{assumption}

The following key identity connects the auxiliary problem~\eqref{eq:sat-a}--\eqref{eq:sat-b} to the residual $\br_h$.

\begin{lemma}
\label{lem:key-identity}
For all $K\in\CT_h$,
\begin{equation}
\|\nabla(\bz_h-\bw_h)\|_K = \|\nabla\br_h\|_K.
\label{eq:key-identity}
\end{equation}
\end{lemma}
\begin{proof}
Subtracting~\eqref{eq:mixed-a} from~\eqref{eq:sat-a} gives $(\nabla(\bz_h-\bw_h),\nabla\bv_K)_K=(\nabla\br_h,\nabla\bv_K)_K$ for all $\bv_K\in(I-P_h)V^{l+1}_h\big|_K$. By~\eqref{eq:sat-b} and~\eqref{eq:mixed-c}, $P_h(\bz_h-\bw_h)=\bzero$, so both $\bz_h-\bw_h$ and $\br_h$ belong to $(I-P_h)V^{l+1}_h\big|_K$, on which $\|\nabla(\cdot)\|_K$ is a norm, since $V_h$ contains the piecewise constants. Hence $\bz_h-\bw_h=\br_h$ on $K$, which gives~\eqref{eq:key-identity}.
\end{proof}

Combined with Assumption~\ref{ass:saturation} and the triangle inequality, Lemma~\ref{lem:key-identity} yields
\begin{equation}
\|\nabla(\bu-\bw_h)\|_0 \leq \frac{1}{1-\delta}\|\nabla\br_h\|_0.
\label{eq:saturated-bound}
\end{equation}

Estimate~\eqref{eq:saturated-bound} controls the full broken gradient of the error, that is, the volume part of $\|\bu-\bw_h\|_{h,\nabla}$, and is the starting point for the compressible-regime reliability estimate established next.

In the analysis below, we also use the following bounds on the elasticity and compliance tensors. Since $\mathcal{A}\btau=2\mu\btau+\lambda\,\text{tr}(\btau)I$, the triangle inequality, together with $\|\text{tr}(\btau)I\|_0\le d\|\btau\|_0$, gives
\begin{equation}\label{eq:A_bound}
    \|\mathcal{A}\btau\|_0 \leq (2\mu+d\lambda)\|\btau\|_0 \quad \forall\btau\in [L^2(\Omega)]^{d\times d}.
\end{equation}
Similarly, $\CC\btau=\frac{1}{2\mu}\big(\btau-\frac{\lambda}{2\mu+d\lambda}\text{tr}(\btau)I\big)$, the same bound on $\|\text{tr}(\btau)I\|_0$, and $\frac{2\mu+2d\lambda}{2\mu+d\lambda}\le2$ give
\begin{equation}\label{eq:C_bound}
    \|\CC\btau\|_0 \leq \mu^{-1}\|\btau\|_0 \quad \forall\btau\in [L^2(\Omega)]^{d\times d}.
\end{equation}

\subsection{Reliability and efficiency estimates}

We are now in a position to state the reliability and efficiency estimates for the indicator $\eta_h$.

\begin{theorem}[Reliability, compressible regime]
\label{thm:reliability-compressible}
Let $(\bsigma,\brho,\bu)$ solve~\eqref{eq:dual-variational}; let $(\bsigma_h,\brho_h,\bu_h)$ solve~\eqref{eq:3.3}; and let $(\br_h,\bw_h)$ solve the residual minimization problem of Section~\ref{sec:residual-min}. If Assumption~\ref{ass:saturation} holds, then
\begin{equation}
    \mu^{-1/2}\|\bsigma-\bsigma_h\|_0 + \mu^{1/2}\|\brho-\brho_h\|_0 + \mu^{1/2}\|\bu-\bw_h\|_{h,\nabla} \le C(\mu,\lambda,\delta)\,\eta_h,
    \label{eq:reliability-compressible}
\end{equation}
where
\begin{equation}
C(\mu,\lambda,\delta)\eq\frac{1}{1-\delta}\Bigg[(3-\delta)+\frac{(2\mu+d\lambda)(2-\delta)}{\mu}\Bigg]+\frac{3}{4}.
\label{eq:explicit-constant}
\end{equation}
\end{theorem}
\begin{proof}
Since $\nabla\bu=\CC\bsigma+\brho$ by~\eqref{eq:2.15}, we have, elementwise,
\begin{equation*}
\nabla(\bu-\bw_h) - \big(\CC\bsigma_h+\brho_h-\nabla\bw_h\big) = \CC(\bsigma-\bsigma_h)+(\brho-\brho_h).
\end{equation*}
Since $\bsigma$ is symmetric, $\brho$ and $\brho_h$ are skew-symmetric, and, $\CC\btau-\frac{1}{2\mu}\btau$ being a multiple of $I$ that depends only on $\tr(\btau)=\tr(\sym(\btau))$, $\sym(\CC\btau)=\CC\sym(\btau)$ and $\skw(\CC\btau)=\frac{1}{2\mu}\skw(\btau)$ for all $\btau\in[L^2(\Omega)]^{d\times d}$, the symmetric and skew-symmetric parts of the right-hand side are $\CC(\bsigma-\sym(\bsigma_h))$ and $(\brho-\brho_h)-\frac{1}{2\mu}\skw(\bsigma_h)=(\brho-\brho_h)-\frac{1}{4\mu}(\bsigma_h-\bsigma_h^{\T})$, respectively. By~\eqref{eq:sym-skw-bound} and the triangle inequality,
\begin{align*}
\|\CC(\bsigma-\sym(\bsigma_h))\|_0 &\le \|\nabla(\bu-\bw_h)\|_0 + \|\CC\bsigma_h+\brho_h-\nabla\bw_h\|_0,\\
\|\brho-\brho_h\|_0 &\le \|\nabla(\bu-\bw_h)\|_0 + \|\CC\bsigma_h+\brho_h-\nabla\bw_h\|_0 + \tfrac{1}{4\mu}\|\bsigma_h-\bsigma_h^{\T}\|_0.
\end{align*}
Since $\bsigma-\sym(\bsigma_h)=\CA\CC(\bsigma-\sym(\bsigma_h))$, applying~\eqref{eq:A_bound} with $\btau=\CC(\bsigma-\sym(\bsigma_h))$, and writing $\bsigma-\bsigma_h=(\bsigma-\sym(\bsigma_h))-\frac{1}{2}(\bsigma_h-\bsigma_h^{\T})$, we also get
\begin{equation*}
\|\bsigma-\bsigma_h\|_0 \le (2\mu+d\lambda)\big(\|\nabla(\bu-\bw_h)\|_0 + \|\CC\bsigma_h+\brho_h-\nabla\bw_h\|_0\big) + \tfrac{1}{2}\|\bsigma_h-\bsigma_h^{\T}\|_0.
\end{equation*}
By~\eqref{eq:saturated-bound}, and since each of $\mu^{1/2}\|\nabla\br_h\|_0$, $\mu^{1/2}\|\CC\bsigma_h+\brho_h-\nabla\bw_h\|_0$, and $\mu^{-1/2}\|\bsigma_h-\bsigma_h^{\T}\|_0$ is bounded by $\eta_h$, it follows that
\begin{equation*}
\mu^{-1/2}\|\bsigma-\bsigma_h\|_0\le\Bigg[\frac{(2\mu+d\lambda)(2-\delta)}{\mu(1-\delta)}+\frac12\Bigg]\eta_h,\quad
\mu^{1/2}\|\brho-\brho_h\|_0\le\Bigg[\frac{2-\delta}{1-\delta}+\frac14\Bigg]\eta_h.
\end{equation*}
Finally, since $\bu\in[H^1_D(\Omega)]^d$ vanishes on $\Gamma_D$ and is single-valued across interior edges/faces, the jump and Dirichlet contributions to $\|\bu-\bw_h\|_{h,\nabla}$ coincide exactly with those of $\bw_h$, so that, by~\eqref{eq:saturated-bound} and since $(1-\delta)^{-2}\ge1$,
\begin{align*}
    \mu\|\bu-\bw_h\|_{h,\nabla}^2&\le\frac{\mu}{(1-\delta)^2}\Bigg(\|\nabla\br_h\|_0^2+\sum_{E\in\CE_h^0}h_E^{-1}\|\jmp{\bw_h}\|_E^2+\sum_{E\in\CE_h^D}h_E^{-1}\|\bw_h\|_E^2\Bigg)\\
    &\le\frac{1}{(1-\delta)^2}\eta_h^2,    
\end{align*}
where the last step uses that these are separate contributions to $\eta_h^2$ in~\eqref{eq:etaK}, each interior edge or face being shared by two elements, with weight $\mu/2$ from each. Adding the three bounds gives~\eqref{eq:reliability-compressible}.
\end{proof}

Since $\|\cdot\|_h\le\|\cdot\|_{h,\nabla}$, estimate~\eqref{eq:reliability-compressible}, and likewise~\eqref{eq:reliability-incompressible}, also holds with $\|\bu-\bw_h\|_h$ in place of $\|\bu-\bw_h\|_{h,\nabla}$.

The reliability estimate of Theorem~\ref{thm:reliability-compressible} degenerates as $\lambda\to\infty$: the constant~\eqref{eq:explicit-constant} grows like $(2\mu+d\lambda)/\mu$, as a direct consequence of~\eqref{eq:A_bound}, reflecting the loss of control over the trace of the stress by the compliance tensor. In plane stress, since $\lambda^*<2\mu$, this constant remains bounded, and Theorem~\ref{thm:reliability-compressible} already suffices. We next derive an estimate with a constant that remains bounded as $\lambda\to\infty$.

Such an estimate relies on the primal-mixed variational formulation of~\eqref{eq:full_system}, equivalent to~\eqref{eq:dual-variational}: find $(\bsigma,\brho,\bu)\in\BE(\Omega)\eq[L^2(\Omega)]^{d\times d}\times[L^2(\Omega)]^{d\times d}_{\skw}\times[H^1_D(\Omega)]^d$ such that
\begin{equation}\label{eq:2.21}
\CB(\bsigma,\brho,\bu;\btau,\bzeta,\bv) = (\bff,\bv) + \langle \bg,\bv\rangle_{\Gamma_N} \qquad \forall(\btau,\bzeta,\bv)\in\BE(\Omega),
\end{equation}
with the bilinear form
\begin{subequations}
\begin{align}
\CB(\bsigma,\brho,\bu;\btau,\bzeta,\bv) &\eq (\CC\bsigma,\btau) + \tilde b(\btau;\brho,\bu) - \tilde b(\bsigma;\bzeta,\bv), \label{eq:bilinear}\\
\tilde b(\btau;\bzeta,\bv) &\eq (\btau,\bzeta)-(\btau,\nabla\bv). \label{eq:btilde}
\end{align}
\end{subequations}
In the norm
\begin{equation}
|\!|\!|(\btau,\bzeta,\bv)|\!|\!|^2 \eq \mu^{-1}\|\btau\|_0^2 + \mu\|\varepsilon(\bv)\|_0^2 + \mu\|\bzeta\|_0^2,
\label{eq:robust-norm}
\end{equation}
the following robust stability estimate holds.

\begin{theorem}[Robust stability]
\label{thm:robust-stability}
It holds that
\begin{equation}
\sup_{(\btau,\bzeta,\bv)\in\BE(\Omega)} \frac{\CB(\bsigma,\brho,\bu;\btau,\bzeta,\bv)}{|\!|\!|(\btau,\bzeta,\bv)|\!|\!|} \;\gtrsim\; |\!|\!|(\bsigma,\brho,\bu)|\!|\!| \qquad \forall(\bsigma,\brho,\bu)\in\BE(\Omega),
\label{eq:robust-stability}
\end{equation}
with a constant that remains bounded as $\lambda\to\infty$.
\end{theorem}
\begin{proof}
See \cite[Theorem~2]{LS2024}.
\end{proof}

To apply Theorem~\ref{thm:robust-stability}, we require a continuous, $H^1_D(\Omega)$-conforming counterpart of the postprocessed displacement $\bw_h$. To this end, let $\osw:V_h^l\to V_h^l\cap[H^1_D(\Omega)]^d$ denote the Oswald averaging operator, which satisfies the approximation property
\begin{equation}\label{eq:oswald_approx}
    \Bigg(\sum_{K\in\CT_h}h_K^{-2}\|\bv-\osw(\bv)\|_K^2\Bigg)^{1/2} \lesssim |\bv|_J \quad \forall\bv\in V_h^l,
\end{equation}
where $|\bv|_J\eq\big(\sum_{E\in\CE_h^0}h_E^{-1}\|\jmp{\bv}\|_E^2+\sum_{E\in\CE_h^D}h_E^{-1}\|\bv\|_E^2\big)^{1/2}$; see~\cite[Theorem 2.1]{KarakashianPascal2007}. Note that $\mu|\bw_h|_J^2$ is precisely the sum of the jump and Dirichlet contributions to $\eta_h^2$ in~\eqref{eq:etaK}.

For the oscillation terms we use the standard quantities
\begin{equation}
\mathrm{osc}(\bff)\eq\Big(\sum_{K\in\CT_h}h_K^2\|\bff-P_h\bff\|_K^2\Big)^{1/2}, \qquad
\mathrm{osc}(\bg)\eq\Big(\sum_{E\in\CE_h^N}h_E\|\bg-Q_E\bg\|_E^2\Big)^{1/2}.
\label{eq:osc}
\end{equation}

\begin{theorem}[Reliability, incompressible limit]
\label{thm:reliability-incompressible}

Let $(\bsigma,\brho,\bu)$ solve~\eqref{eq:dual-variational}; let $(\bsigma_h,\brho_h,\bu_h)$ solve~\eqref{eq:3.3}; and let $(\br_h,\bw_h)$ solve the residual minimization problem of Section~\ref{sec:residual-min}. Then, for every $\lambda>0$, it holds that

\begin{align}
\mu^{-1/2}\|\bsigma-\bsigma_h\|_0 + \mu^{1/2}\|\brho-\brho_h\|_0\,+&\,\mu^{1/2}\|\bu-\bw_h\|_{h,\nabla}\nonumber\\
&\lesssim \eta_h + \mu^{-1/2}\big(\mathrm{osc}(\bff)+\mathrm{osc}(\bg)\big),
\label{eq:reliability-incompressible}
\end{align}
with a constant that remains bounded as $\lambda\to\infty$.
\end{theorem}
\begin{proof}
The argument follows that of~\cite[Theorem~8]{LS2024}, with the averaged postprocessed displacement used there replaced by $\osw\bw_h$: Theorem~\ref{thm:robust-stability} is applied to the error triple $(\bsigma-\bsigma_h,\brho-\brho_h,\bu-\osw\bw_h)\in\BE(\Omega)$. Two terms require separate treatment. The first, $\CC\bsigma_h+\brho_h-\nabla(\osw\bw_h)$, is controlled by the triangle inequality and the approximation property of $\osw$: writing $\CC\bsigma_h+\brho_h-\nabla(\osw\bw_h)=\big(\CC\bsigma_h+\brho_h-\nabla\bw_h\big)+\nabla(\bw_h-\osw\bw_h)$, the first term is a component of $\eta_K$, while the second is bounded by~\eqref{eq:oswald_approx} and an inverse inequality applied to the polynomial $\bw_h-\osw(\bw_h)$ on each element. The second, the coupling term $\frac{1}{2}(\bsigma_h-\bsigma_h^{\T},\bzeta)$ arising from the skew-symmetric part of the bilinear form, cannot be treated via averaging, since it is not an $H^1$-type quantity; instead, by Cauchy--Schwarz and $|\!|\!|(\btau,\bzeta,\bv)|\!|\!|\ge\mu^{1/2}\|\bzeta\|_0$, it is bounded directly by
\begin{equation*}
\frac{1}{2}\frac{|(\bsigma_h-\bsigma_h^{\T},\bzeta)|}{|\!|\!|(\btau,\bzeta,\bv)|\!|\!|} \le \frac{\|\bsigma_h-\bsigma_h^{\T}\|_0\,\|\bzeta\|_0}{2\mu^{1/2}\|\bzeta\|_0} = \frac{\mu^{-1/2}}{2}\|\bsigma_h-\bsigma_h^{\T}\|_0,
\end{equation*}
which is bounded by the asymmetry defect in~\eqref{eq:etaK}. Theorem~\ref{thm:robust-stability} thus controls $\|\varepsilon(\bu-\osw\bw_h)\|_0$, and, by~\eqref{eq:oswald_approx} and an inverse inequality, $\|\nabla(\bw_h-\osw\bw_h)\|_0\lesssim|\bw_h|_J$, so that $\|\varepsilon(\bu-\bw_h)\|_0$ is controlled as well. The skew-symmetric part of $\nabla(\bu-\bw_h)$ is recovered from $\omega(\bu)=\brho$ and $\skw(\CC\bsigma_h)=\frac{1}{2\mu}\skw(\bsigma_h)$, which give
\begin{equation*}
\omega(\bu-\bw_h)=(\brho-\brho_h)+\skw\big(\CC\bsigma_h+\brho_h-\nabla\bw_h\big)-\tfrac{1}{4\mu}(\bsigma_h-\bsigma_h^{\T}),
\end{equation*}
so that $\|\omega(\bu-\bw_h)\|_0$ is bounded by quantities already controlled. The jump and Dirichlet contributions to $\|\bu-\bw_h\|_{h,\nabla}$ coincide exactly with those of $\bw_h$, that is, with $|\bw_h|_J$, as in the proof of Theorem~\ref{thm:reliability-compressible}. The oscillation terms are handled as in~\cite[(6.6)--(6.9)]{LS2024}.
\end{proof}

\begin{remark}[Treatment of the asymmetry defect]
In the proof of Theorem~\ref{thm:reliability-incompressible}, the asymmetry of $\bsigma_h$ enters only through the coupling term $\frac12(\bsigma_h-\bsigma_h^{\T},\bzeta)$, which, thanks to the rotation component of the norm~\eqref{eq:robust-norm}, is bounded directly by the asymmetry defect in~\eqref{eq:etaK}, with no additional inequality required.
\end{remark}

\begin{remark}[A single estimator for both regimes]
Theorems~\ref{thm:reliability-compressible} and~\ref{thm:reliability-incompressible} bound the same error quantity by the same indicator $\eta_h$, up to data oscillation in the latter: the former yields the explicit constant~\eqref{eq:explicit-constant}, which grows with $\lambda$, and the latter a constant that remains bounded as $\lambda\to\infty$. The hypercircle estimator of~\cite[Theorem~7]{LS2024} is formulated in the norm induced by $\CC$, which degenerates in the incompressible limit; for that limit, \cite{LS2024} uses instead the estimator of~\cite[Theorem~8]{LS2024}, uniformly valid in $\lambda$, with which ours shares the asymmetry term and a stress residual of the same form.
\end{remark}

We conclude with a local efficiency estimate, showing that the indicator $\eta_K$ is bounded above by the local error, with an explicit constant independent of $h$, $\mu$, and $\lambda$, and without any saturation assumption.

\begin{theorem}[Local efficiency]
\label{thm:efficiency}
Let $(\bsigma,\brho,\bu)$ solve~\eqref{eq:dual-variational}; let $(\bsigma_h,\brho_h,\bu_h)$ solve~\eqref{eq:3.3}; and let $(\br_h,\bw_h)$ solve the residual minimization problem of Section~\ref{sec:residual-min}. Then, for all $K\in\CT_h$,
\begin{equation}
\eta_K \le \sqrt{10}\,\Big(\mu^{-1}\|\bsigma-\bsigma_h\|_K^2 + \mu\|\brho-\brho_h\|_K^2 + \mu\|\bu-\bw_h\|_{h,\nabla,K}^2\Big)^{1/2}.
\label{eq:efficiency}
\end{equation}
\end{theorem}
\begin{proof}
Let $G_K\eq(\CC\bsigma_h+\brho_h-\nabla\bw_h)|_K$. Since $\CC\bsigma+\brho=\nabla\bu$,
\begin{equation*}
G_K=-\CC(\bsigma-\bsigma_h)-(\brho-\brho_h)+\nabla(\bu-\bw_h),
\end{equation*}
so that, by~\eqref{eq:C_bound},
\begin{equation*}
\mu^{1/2}\|G_K\|_K\le R_K\eq\mu^{-1/2}\|\bsigma-\bsigma_h\|_K+\mu^{1/2}\|\brho-\brho_h\|_K+\mu^{1/2}\|\nabla(\bu-\bw_h)\|_K.
\end{equation*}
Testing~\eqref{eq:mixed-a} with $\bv_K=\br_h|_K$ gives $\|\nabla\br_h\|_K^2=(G_K,\nabla\br_h)_K\le\|G_K\|_K\|\nabla\br_h\|_K$, whence $\mu^{1/2}\|\nabla\br_h\|_K\le R_K$ as well. Moreover, since $\bsigma$ is symmetric, $\bsigma_h-\bsigma_h^{\T}=(\bsigma_h-\bsigma)-(\bsigma_h-\bsigma)^{\T}$, so that $\mu^{-1/2}\|\bsigma_h-\bsigma_h^{\T}\|_K\le2\mu^{-1/2}\|\bsigma-\bsigma_h\|_K$. Finally, as in the proof of Theorem~\ref{thm:reliability-compressible}, the jump and Dirichlet contributions to $\eta_K$ coincide exactly with those of $\mu\|\bu-\bw_h\|_{h,\nabla,K}^2$. Using $R_K^2\le3\big(\mu^{-1}\|\bsigma-\bsigma_h\|_K^2+\mu\|\brho-\brho_h\|_K^2+\mu\|\nabla(\bu-\bw_h)\|_K^2\big)$ and collecting terms gives~\eqref{eq:efficiency}.
\end{proof}

The constant in Theorem~\ref{thm:efficiency} is uniform with respect to $h$, $\mu$, and $\lambda$; hence the estimate holds in both the compressible regime and the incompressible limit, in contrast with the reliability estimate of Theorem~\ref{thm:reliability-compressible}, whose constant grows with $\lambda$.

\begin{remark}[Scope of the analysis]\label{rem:scope}
The a posteriori analysis above uses the discrete problem~\eqref{eq:3.3} only through three properties: the inclusion $\ddiv S_h\subset V_h$, which yields $\ddiv\bsigma_h+P_h\bff=\bzero$, and the traction condition $\bsigma_h\bn=Q_h\bg$ on $\Gamma_N$, both entering the oscillation terms of Theorem~\ref{thm:reliability-incompressible}; and the inclusion of the piecewise constants in $V_h$, used in Lemma~\ref{lem:key-identity}. Proposition~\ref{prop:equivalence} and Theorems~\ref{thm:reliability-compressible}, \ref{thm:reliability-incompressible}, and~\ref{thm:efficiency} therefore remain valid for other stable weakly symmetric discretizations with these properties, such as PEERS~\cite{ArnoldBrezziDouglas1984} and the other simplicial elements analyzed in~\cite{Lee2016}. Only the a priori rates of Theorem~\ref{thm:postproc-apriori} rely on the family-specific estimates of~\cite{LS2024}.
\end{remark}

\section{Numerical examples}
\label{sec:numerical_results}

We validate the theoretical results of the previous sections with three two-dimensional numerical examples, implemented in Firedrake~\cite{FiredrakeUserManual}. Two weakly symmetric families are considered: the AFW element, used throughout as the primary family, and the SGG element, used in the smooth example and in the L-shaped domain example to verify the sharper convergence rate predicted for this family alone in Remark~\ref{rem:superconvergence}. No comparison between families is intended; each is reported on its own terms.

For each example, we report the relative errors
\begin{equation}
e_0^{\bsigma} = \frac{\|\bsigma-\bsigma_h\|_0}{\|\bsigma\|_0}, \qquad
e_0^{\brho} = \frac{\|\brho-\brho_h\|_0}{\|\brho\|_0}, \qquad
e_h^{\bu} = \frac{\|\bu-\bw_h\|_{h,\nabla}}{\|\bu\|_{h,\nabla}},
\end{equation}
and the effectivity index
\begin{equation}\label{eq:effectivity_index}
    c_{\text{eff}} = \frac{\eta_h}{e_{\rm total}},\qquad e_{\rm total}\eq\big(\mu^{-1}\|\bsigma-\bsigma_h\|_0^2+\mu\|\brho-\brho_h\|_0^2+\mu\|\bu-\bw_h\|_{h,\nabla}^2\big)^{1/2}.
\end{equation}
Adaptive computations follow the standard loop
\[
\text{SOLVE} \to \text{ESTIMATE} \to \text{MARK} \to \text{REFINE} \to \text{SOLVE} \to \cdots,
\]
with the maximum marking strategy, as in~\cite{LS2023,LS2024}: recalling the local indicator $\eta_K$ of~\eqref{eq:etaK}, an element $K\in\CT_h$ is refined whenever $\eta_K^2\ge\theta\max_{K'\in\CT_h}\eta_{K'}^2$, $\theta=1/4$. We denote by $N$ the total number of degrees of freedom of the discrete problem~\eqref{eq:3.3}. We specify the material by its Young's modulus $E$ and Poisson's ratio $\nu\in(0,1/2)$, related to the Lam\'e parameters by
\begin{equation}
\lambda = \frac{E\nu}{(1+\nu)(1-2\nu)}, \qquad \mu = \frac{E}{2(1+\nu)}.
\end{equation}
All examples are posed in plane strain, where the incompressible limit $\lambda\to\infty$ corresponds to $\nu\to1/2$.

\subsection{A smooth solution}
\label{sec:num-smooth}

To verify the predicted convergence rates, we let $\Omega=(0,1)^2$, $\Gamma_D=\partial\Omega$, $\Gamma_N=\emptyset$, and take the exact displacement
\begin{align*}
u_1(x,y) &= \pi\sin^2(\pi x)\sin(2\pi y)+\frac{1}{\lambda}\sin(\pi x)\sin(\pi y),\\
u_2(x,y) &= -\pi\sin(2\pi x)\sin^2(\pi y)+\frac{1}{\lambda}\sin(\pi x)\sin(\pi y),
\end{align*}
with body force $\bff=-\ddiv\bsigma$, stress $\bsigma=\CC^{-1}\varepsilon(\bu)$, and exact rotation $\brho=\omega(\bu)$. Since $\lambda\ddiv\bu=\pi\sin(\pi(x+y))$ is independent of $\lambda$, both $\bsigma$ and $\bff$ remain bounded as $\lambda\to\infty$. We fix $E=1$.

As noted following Theorem~\ref{thm:apriori}, for a sufficiently smooth exact solution the expected rate for the stress in the $L^2$-norm is $\CO(h^k)$ for AFW; by Corollary~\ref{cor:apriori-residual}, the postprocessed displacement $\bw_h$ attains the corresponding order in the broken norm $\|\cdot\|_{h,\nabla}$. Table~\ref{tab:smooth-afw} confirms that the expected rate is attained for both values of $\nu$, and that the rate is attained without degradation as $\nu\to1/2$, confirming that the method is free of the locking phenomenon that affects standard primal discretizations in the incompressible limit.
\begin{table}[h]
\centering
\small
\caption{Smooth example: relative errors and orders of convergence (with respect to $h$) for the AFW method for $k=4$, uniform refinement, $\nu=0.3$ (top) and $\nu=0.49999$ (bottom).}
\label{tab:smooth-afw}
\begin{tabular}{lccccccccc}
\hline
$N$ & $e_0^{\bsigma}$ & (oc) & $e_0^{\brho}$ & (oc) & $e_h^{\bu}$ & (oc) & $\eta_h$ & (oc) & $c_{\rm eff}$\\
\hline
5520   & 4.59e-4 & $-$  & 8.23e-4 & $-$  & 4.14e-4 & $-$  & 9.51e-3 & $-$ &  0.647 \\
21840  & 2.85e-5 & 4.01 & 5.34e-5 & 3.94 & 2.67e-5 & 3.95 & 5.90e-4 & 4.01 & 0.632 \\
86880  & 1.77e-6 & 4.01 & 3.38e-6 & 3.98 & 1.69e-6 & 3.98 & 3.65e-5 & 4.02 & 0.622 \\ 
195120 & 3.49e-7 & 4.01 & 6.70e-7 & 3.99 & 3.35e-7 & 3.99 & 7.17e-6 & 4.01 & 0.620 \\
346560 & 1.10e-7 & 4.00 & 2.12e-7 & 4.00 & 1.06e-7 & 4.00 & 2.27e-6 & 4.01 & 0.619 \\
\hline
5520   & 5.43e-4 & $-$  & 8.51e-4 & $-$  & 4.45e-4 & $-$  & 9.00e-3 & $-$ &  0.649 \\
21840  & 3.41e-5 & 4.00 & 5.50e-5 & 3.95 & 2.85e-5 & 3.96 & 5.61e-4 & 4.00 & 0.636 \\
86880  & 2.12e-6 & 4.00 & 3.48e-6 & 3.98 & 1.80e-6 & 3.98 & 3.47e-5 & 4.01 & 0.627 \\ 
195120 & 4.19e-7 & 4.00 & 6.88e-7 & 3.99 & 3.57e-7 & 3.99 & 6.84e-6 & 4.01 & 0.625 \\
346560 & 1.33e-7 & 4.00 & 2.18e-7 & 4.00 & 1.13e-7 & 4.00 & 2.16e-6 & 4.01 & 0.623 \\
\hline
\end{tabular}
\end{table}

To verify that the SGG family attains the sharper rate predicted for it alone in Remark~\ref{rem:superconvergence}, Table~\ref{tab:smooth-sgg} reports the same quantities for SGG on the same example.

\begin{table}[h]
\centering
\small
\caption{Smooth example: relative errors and orders of convergence (with respect to $h$) for the SGG method for $k=4$, uniform refinement, $\nu=0.3$ (top) and $\nu=0.49999$ (bottom).}
\label{tab:smooth-sgg}
\begin{tabular}{lccccccccc}
\hline
$N$ & $e_0^{\bsigma}$ & (oc) & $e_0^{\brho}$ & (oc) & $e_h^{\bu}$ & (oc) & $\eta_h$ & (oc) & $c_{\rm eff}$\\
\hline
6240   & 9.12e-5 & $-$  & 1.00e-4 & $-$  & 1.12e-4 & $-$  & 6.98e-4 & $-$ &  0.239 \\
24720  & 2.96e-6 & 4.95 & 3.28e-6 & 4.93 & 3.64e-6 & 4.94 & 2.27e-5 & 4.94 & 0.240 \\
98400  & 9.36e-8 & 4.98 & 1.04e-7 & 4.98 & 1.15e-7 & 4.98 & 7.20e-7 & 4.98 & 0.240 \\ 
221040 & 1.24e-8 & 4.99 & 1.38e-8 & 4.99 & 1.52e-8 & 4.99 & 9.52e-8 & 4.99 & 0.240 \\
392640 & 2.94e-9 & 5.00 & 3.28e-9 & 4.99 & 3.62e-9 & 5.00 & 2.26e-8 & 4.99 & 0.240 \\
\hline
6240   & 1.09e-4 & $-$  & 1.04e-4 & $-$  & 1.19e-4 & $-$  & 6.23e-4 & $-$ &  0.227 \\
24720  & 3.51e-6 & 4.95 & 3.40e-6 & 4.93 & 3.87e-6 & 4.94 & 2.02e-5 & 4.95 & 0.226 \\
98400  & 1.11e-7 & 4.99 & 1.08e-7 & 4.98 & 1.22e-7 & 4.98 & 6.37e-7 & 4.98 & 0.226 \\ 
221040 & 1.46e-8 & 5.00 & 1.43e-8 & 4.99 & 1.62e-8 & 4.99 & 8.41e-8 & 4.99 & 0.226 \\
392640 & 3.47e-9 & 5.00 & 3.40e-9 & 4.99 & 3.84e-9 & 5.00 & 2.00e-8 & 5.00 & 0.226 \\
\hline
\end{tabular}
\end{table}
	
Since the exact solution is known, this example also allows us to verify the behaviour of the saturation constant $\delta$ of Assumption~\ref{ass:saturation}. Table~\ref{tab:saturation} reports the computed ratio
\begin{equation}
\delta_h(\nu) \eq \frac{\|\nabla(\bu-\bz_h)\|_0}{\|\nabla(\bu-\bw_h)\|_0}
\end{equation}
for both families, across a range of mesh sizes and Poisson ratios $\nu$. In all cases, $\delta_h(\nu)$ stays below one and is essentially independent of $N$ and $\nu$, consistent with Assumption~\ref{ass:saturation}.

\begin{table}[h]
\centering
\small
\caption{Smooth example: computed saturation ratio $\delta_h(\nu)$, AFW method (top) and SGG method (bottom), for varying $N$ and $\nu$ with $k=4$.}
\label{tab:saturation}
\begin{tabular}{rccccc}
\hline
$N$ & $\nu=0.2$ & $\nu=0.3$ & $\nu=0.4$ & $\nu=0.49$ & $\nu=0.49999$ \\
\hline
5520  & 0.966 & 0.966 & 0.968 & 0.970 & 0.971 \\
21840  & 0.962 & 0.963 & 0.964 & 0.967 & 0.968 \\
86880 & 0.959 & 0.960 & 0.962 & 0.965 & 0.965 \\
195120 & 0.958 & 0.959 & 0.961 & 0.964 & 0.964 \\
346560 & 0.958 & 0.959 & 0.960 & 0.964 & 0.964 \\
\hline
6240  & 0.795 & 0.795 & 0.797 & 0.800 & 0.800 \\
24720  & 0.795 & 0.796 & 0.798 & 0.801 & 0.802 \\
98400 & 0.794 & 0.796 & 0.798 & 0.801 & 0.802 \\
221040 & 0.794 & 0.796 & 0.798 & 0.802 & 0.802 \\
392640 & 0.794 & 0.796 & 0.798 & 0.802 & 0.802 \\
\hline
\end{tabular}
\end{table}

\begin{figure}[h]
\centering
\includegraphics[width=0.49\textwidth]{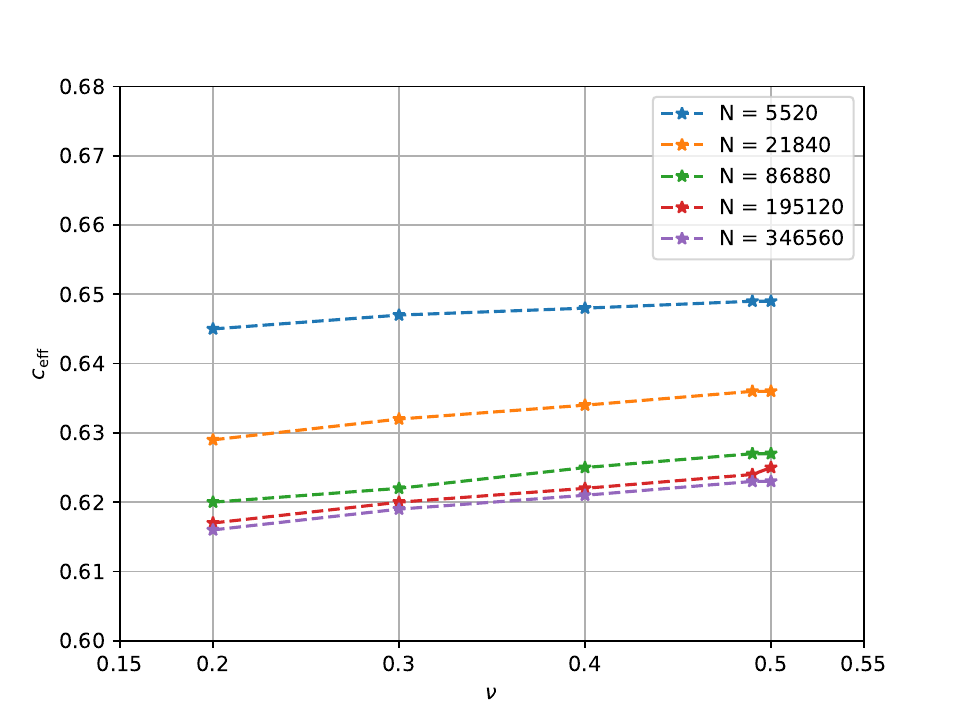}
\includegraphics[width=0.49\textwidth]{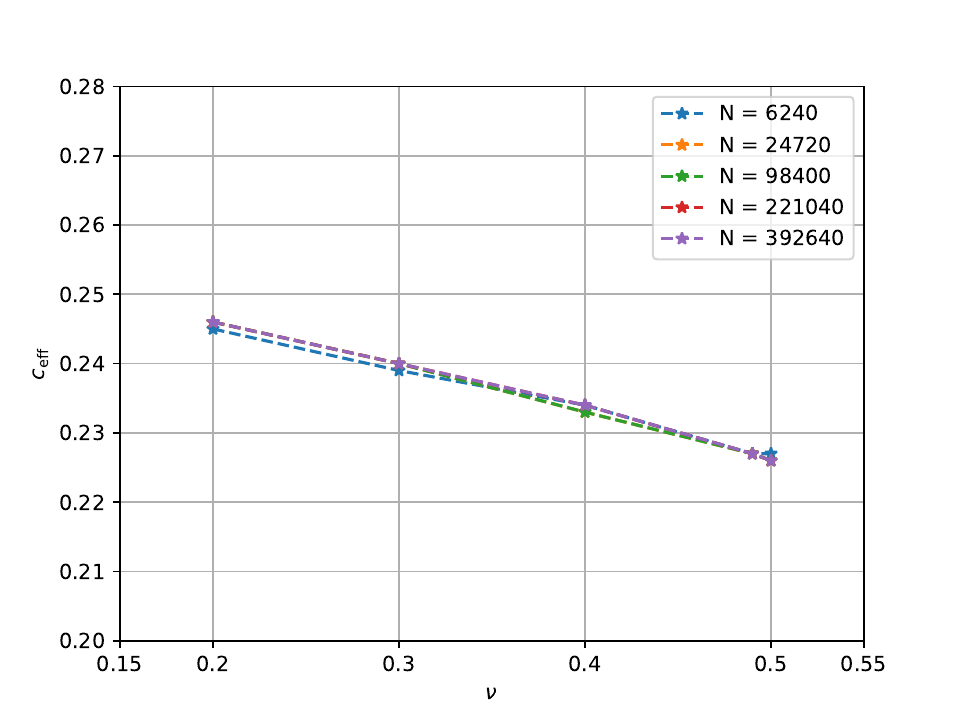}
\caption{Smooth example: effectivity index $c_{\text{eff}}$ as a function of $\nu$, AFW method (left) and SGG method (right) for $k=4$.}
\label{fig:smooth-ceff}
\end{figure}

Figure~\ref{fig:smooth-ceff} shows $c_{\text{eff}}$, as defined in~\eqref{eq:effectivity_index}, as a function of $\nu$, for both element families and for five values ranging from $0.2$ to $0.49999$: for each fixed $N$, the index varies only mildly with $\nu$, and it settles as $N$ increases, in agreement with Theorems~\ref{thm:reliability-incompressible} and~\ref{thm:efficiency}, whose constants remain bounded as $\lambda\to\infty$.

\subsection{An L-shaped domain}
\label{sec:num-Lshape}

Next, we consider the benchmark from~\cite[Section 8.2]{LS2023}. The domain
\[
\Omega = \{(x,y): |x|+|y|\le \sqrt2\,a\} \setminus \{(x,y): |x+a/\sqrt2|+|y|\le a/\sqrt2\},
\]
with $a=1$, has a reentrant corner of interior angle $3\pi/2$, and traction boundary conditions are imposed on the whole boundary, $\Gamma_N=\Gamma$, $\Gamma_D=\emptyset$.  This benchmark, also used in~\cite{LS2024}, lies formally outside the setting $|\Gamma_D|>0$ of Section~\ref{sec:governing_eqns}; as there, rigid body motions are factored out.

Since $\Gamma=\Gamma_N$, the displacement $\bu$ is determined only up to a rigid body motion $\boldsymbol{a}+B\bx$, with $\boldsymbol{a}\in\mathbb{R}^2$ and $B\in\mathbb{R}^{2\times2}$ skew-symmetric, and the rotation $\brho=\omega(\bu)$, correspondingly, only up to $B$; a Lagrange multiplier is introduced in the implementation to fix both. This choice does not affect the indicator: adding a pair $(\boldsymbol{a}+B\bx,B)\in V_h\times R_h$ to $(\bu_h,\brho_h)$ adds $\boldsymbol{a}+B\bx$ to $\bw_h$, by~\eqref{eq:step1a}--\eqref{eq:step1b}, and leaves $\CC\bsigma_h+\brho_h-\nabla\bw_h$, $\br_h$, and the jumps of $\bw_h$ unchanged, so that $\eta_h$ is the same for every such choice. The errors, in contrast, depend on it through $\brho_h$ and $\nabla\bw_h$; we therefore compute them modulo infinitesimal rotations, replacing $(\brho_h,\bw_h)$ by $(\brho_h+B^*,\bw_h+B^*\bx)$, where
\[
B^*\eq\frac{1}{2|\Omega|}\int_\Omega\big((\brho-\brho_h)+\omega(\bu-\bw_h)\big)\,dx
\]
is the constant skew-symmetric matrix minimizing $e_{\rm total}$. In polar coordinates centered at the reentrant corner, the exact displacement is

\begin{align*}
u_x &= \frac{1}{2\mu}r^{\alpha}\big((\kappa-Q(\alpha+1))\cos(\alpha\theta)-\alpha\cos((\alpha-2)\theta)\big),\\
u_y &= \frac{1}{2\mu}r^{\alpha}\big((\kappa+Q(\alpha+1))\sin(\alpha\theta)+\alpha\sin((\alpha-2)\theta)\big),
\end{align*}
with $\alpha=0.544483737$, $Q=0.543075579$, $\kappa=3-4\nu$. We fix $E=1$ and consider $\nu=0.3$ and $\nu=0.49999$, using the AFW and SGG methods with $k=4$. For the study of $c_{\text{eff}}$ below, $\nu$ also takes the values $0.2$, $0.4$ and $0.49$.

\begin{figure}[h]
\centering
 \includegraphics[width=0.49\textwidth]{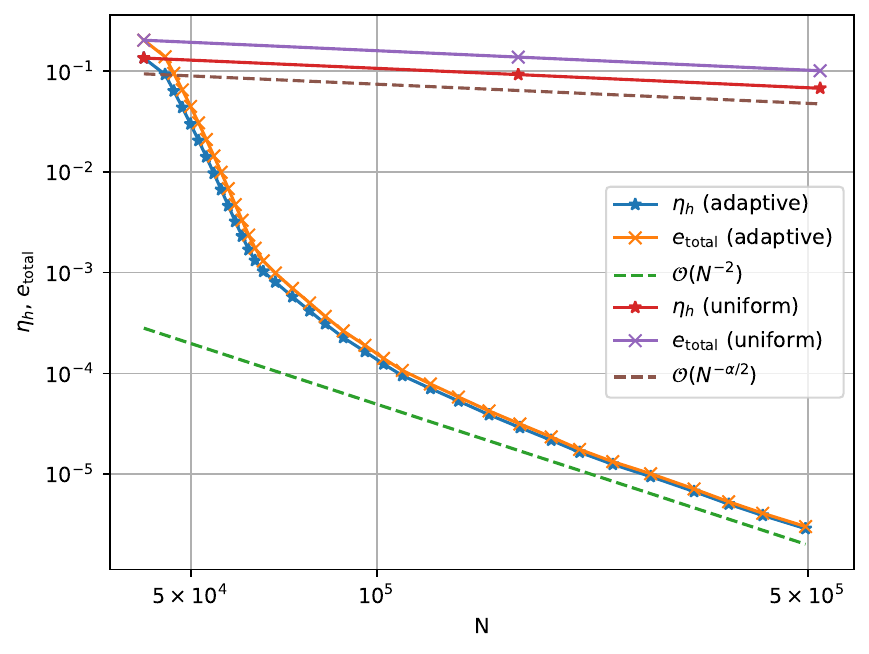}
 \includegraphics[width=0.49\textwidth]{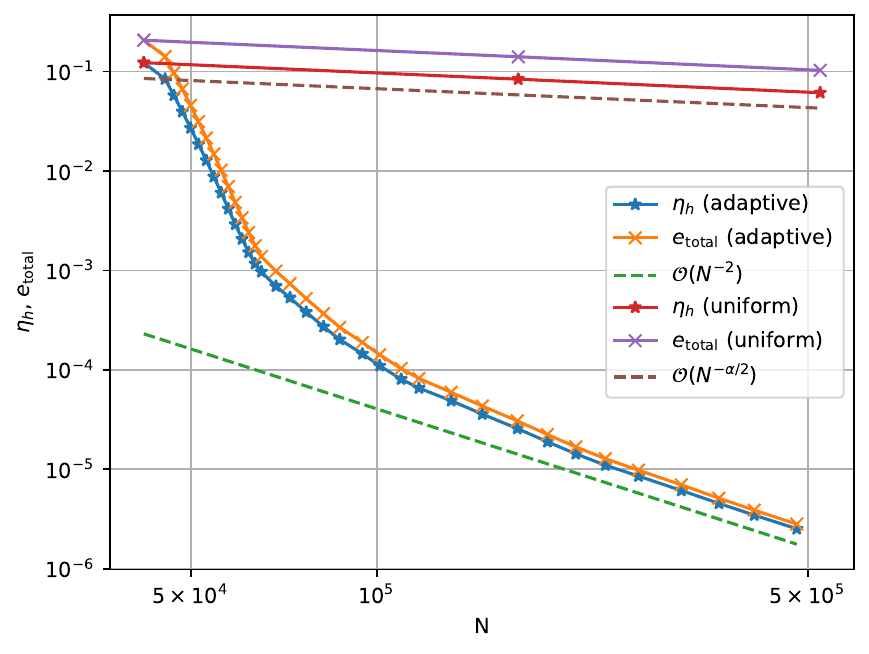}
\caption{L-shaped domain: convergence history of $\eta_h$ and $e_{\rm total}$ under adaptive and uniform refinement, $\nu=0.3$ (left) and $\nu=0.49999$ (right), AFW method, $k=4$.}
\label{fig:Lshape-conv-afw}
\end{figure}

\begin{figure}[h]
\centering
 \includegraphics[width=0.49\textwidth]{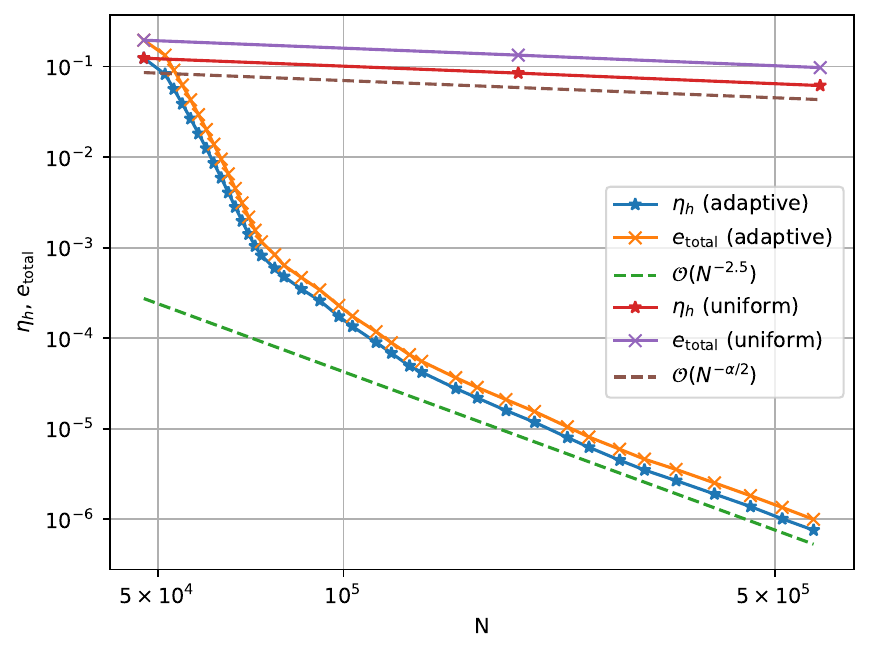}
 \includegraphics[width=0.49\textwidth]{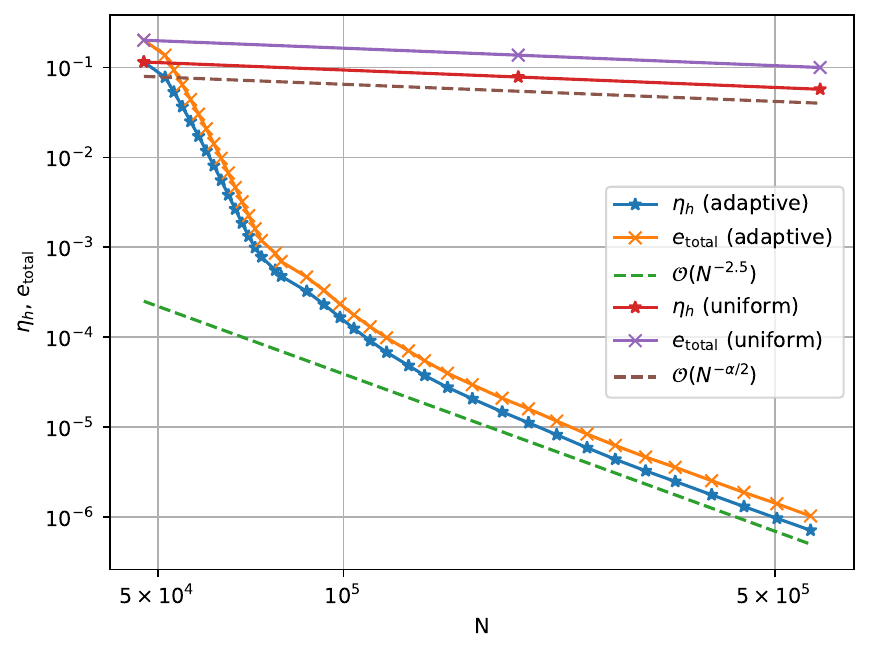}
\caption{L-shaped domain: convergence history of $\eta_h$ and $e_{\rm total}$ under adaptive and uniform refinement, $\nu=0.3$ (left) and $\nu=0.49999$ (right), SGG method, $k=4$.}
\label{fig:Lshape-conv-sgg}
\end{figure}

\begin{figure}[h]
\centering
 \includegraphics[width=\textwidth]{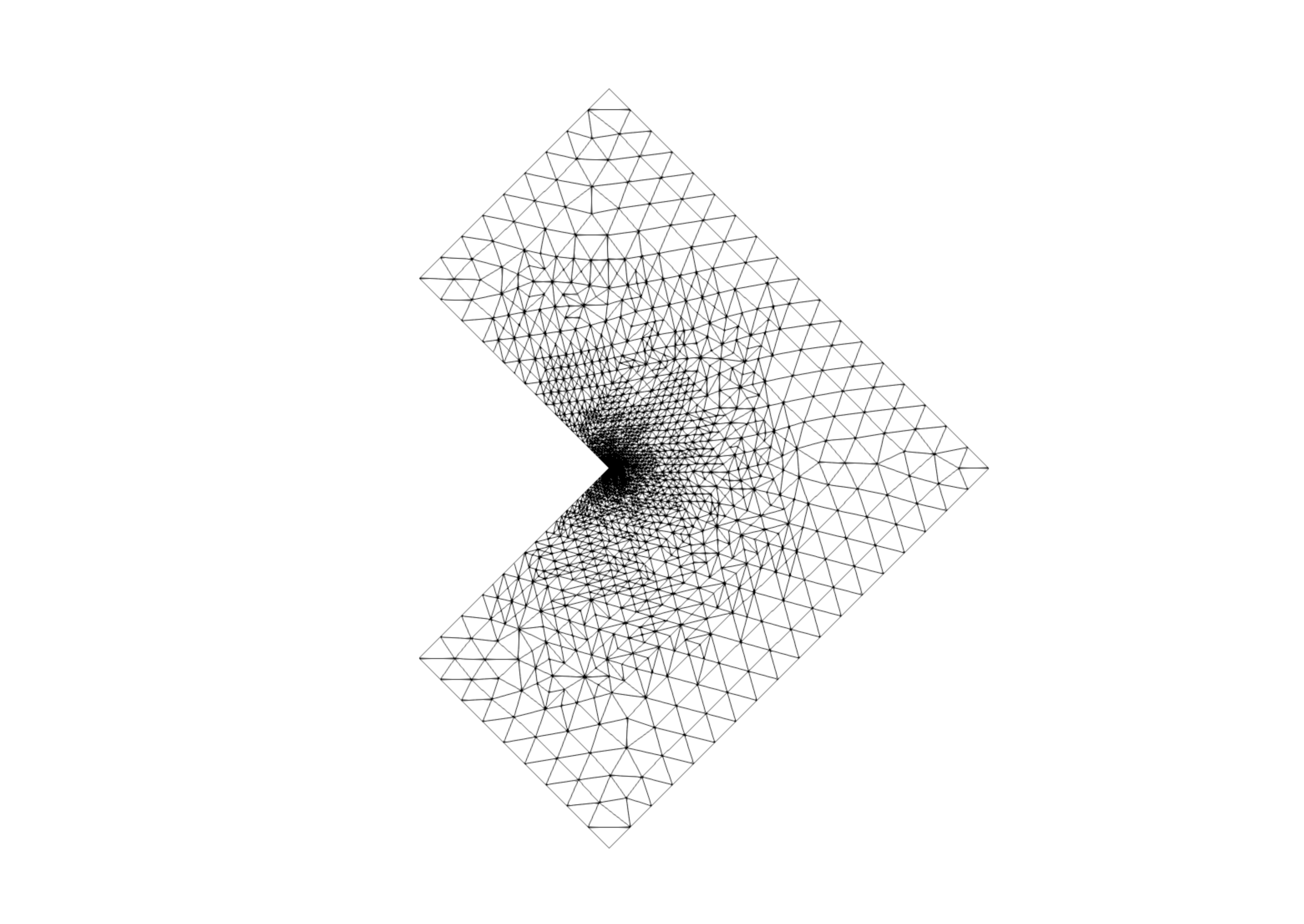}
\caption{L-shaped domain: adaptively refined mesh after 34 refinements, AFW method, $\nu=0.49999$.}
\label{fig:Lshaped-mesh}
\end{figure}

Since $\bu\in H^{1+s}(\Omega)$ only for $s<\alpha$, with $\alpha<0.55$~\cite{SzaboBabuska1991}, uniform refinement yields the suboptimal rate $\CO(N^{-\alpha/2})$, independently of the polynomial degree. Figures~\ref{fig:Lshape-conv-afw} and \ref{fig:Lshape-conv-sgg} compare the convergence histories of $\eta_h$ and $e_{\rm total}$ under adaptive refinement, driven by the indicator of Section~\ref{sec:aposteriori}, and under uniform refinement, for both values of $\nu$ and for both families: adaptivity recovers the optimal rate $\CO(N^{-l/2})$, that is, $\CO(N^{-k/2})$ for AFW and $\CO(N^{-(k+1)/2})$ for SGG, for both quantities. Figure~\ref{fig:Lshaped-mesh} displays a representative adaptively refined mesh, confirming that refinement concentrates at the reentrant corner.

\begin{figure}[h]
\centering
\includegraphics[width=0.49\textwidth]{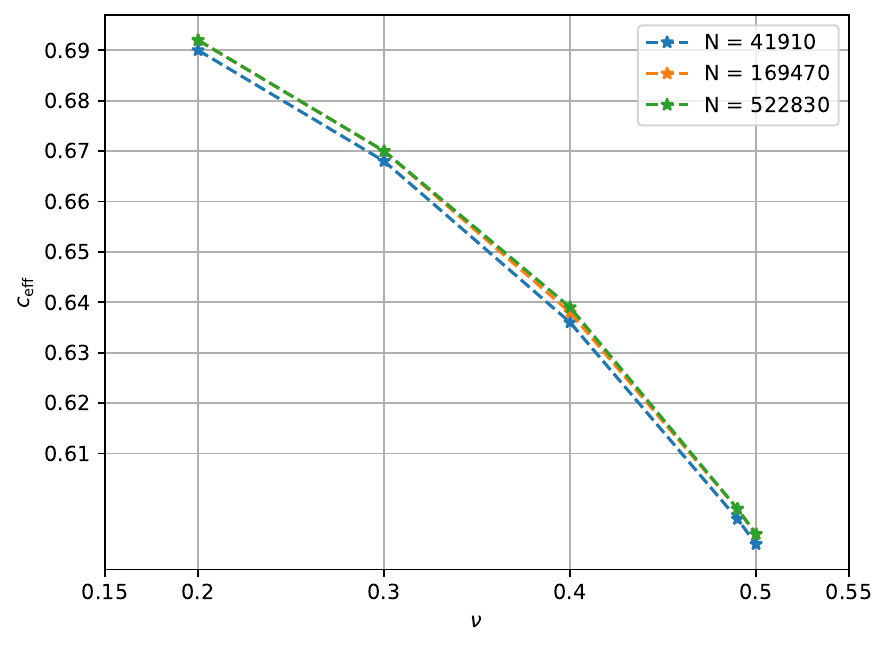}
\includegraphics[width=0.49\textwidth]{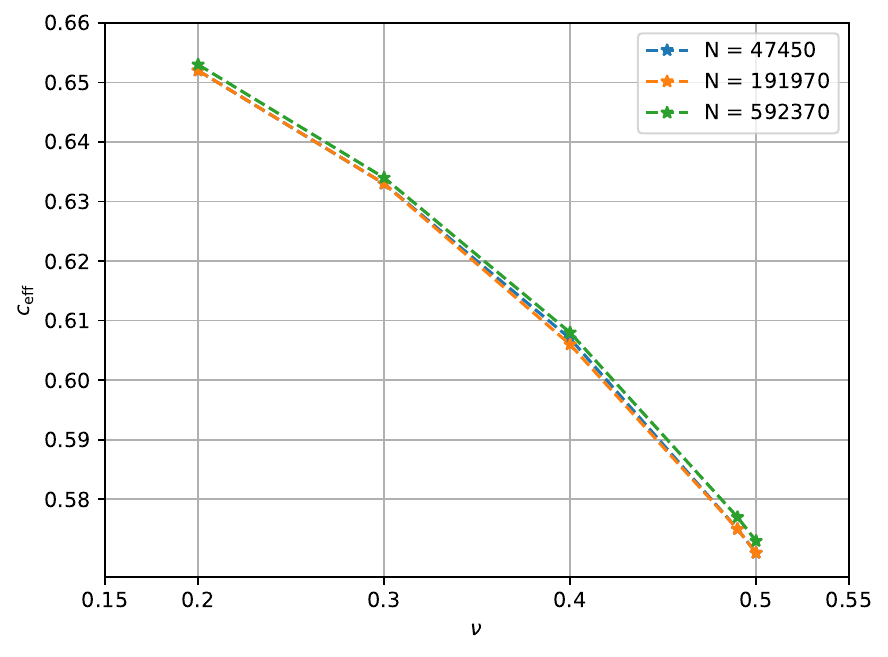}
\caption{L-shaped domain: effectivity index $c_{\text{eff}}$ as a function of $\nu$ on three fixed meshes, AFW method (left) and SGG method (right), $k=4$.}
\label{fig:lshaped-ceff}
\end{figure}

Figure~\ref{fig:lshaped-ceff} shows $c_{\text{eff}}$, as defined in~\eqref{eq:effectivity_index}, as a function of $\nu\in\{0.2,\allowbreak 0.3,\allowbreak 0.4,\allowbreak 0.49,\allowbreak 0.49999\}$, for both families on three fixed meshes, the same for every value of $\nu$: the index decreases moderately as $\nu\to1/2$, by less than $15\%$ over this range, and is nearly independent of $N$, in agreement with Theorems~\ref{thm:reliability-incompressible} and~\ref{thm:efficiency}, whose constants remain bounded as $\lambda\to\infty$.

\subsection{Cook's membrane}
\label{sec:num-cook}

Our last example is Cook's membrane, a standard benchmark in computational mechanics for which no analytical solution is available. The domain is the trapezoid with vertices $(0,0)$, $(48,44)$, $(48,60)$, $(0,44)$, clamped on $\Gamma_D=\text{conv}\{(0,0),(0,44)\}$, subject to a shear traction $\bg=(0,1)$ on $\text{conv}\{(48,44),(48,60)\}$ and traction-free elsewhere, and $\bff=(0,0)$. We take $E=1$ and $\nu\in\{0.3,0.49999\}$, using the AFW method with $k=4$; the incompatible traction data at the corners $(48,44)$ and $(48,60)$ make this a demanding benchmark, especially for $\nu=0.49999$.

\begin{figure}[h!]
\centering
 \includegraphics[width=0.49\textwidth]{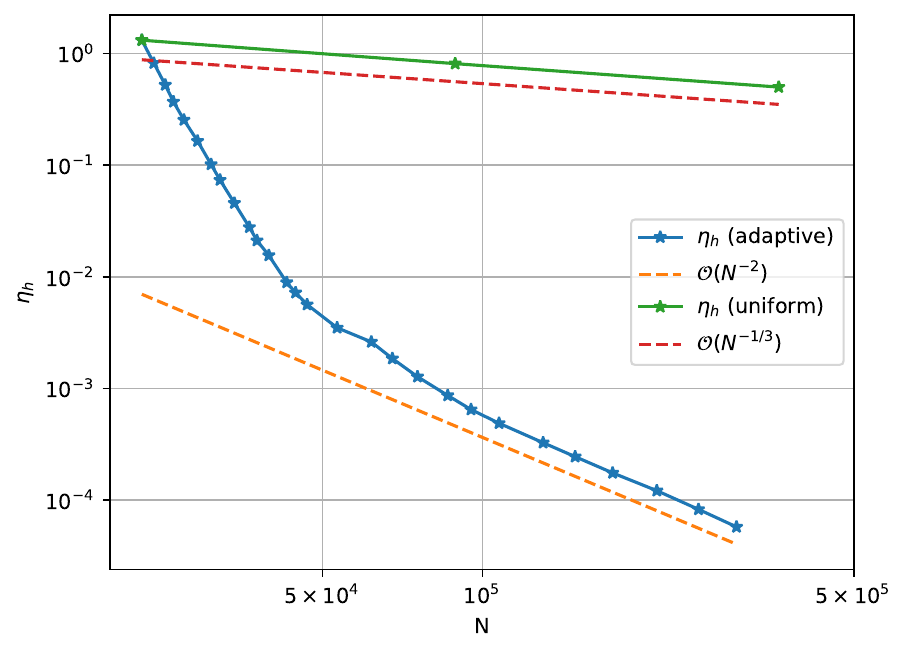}
 \includegraphics[width=0.49\textwidth]{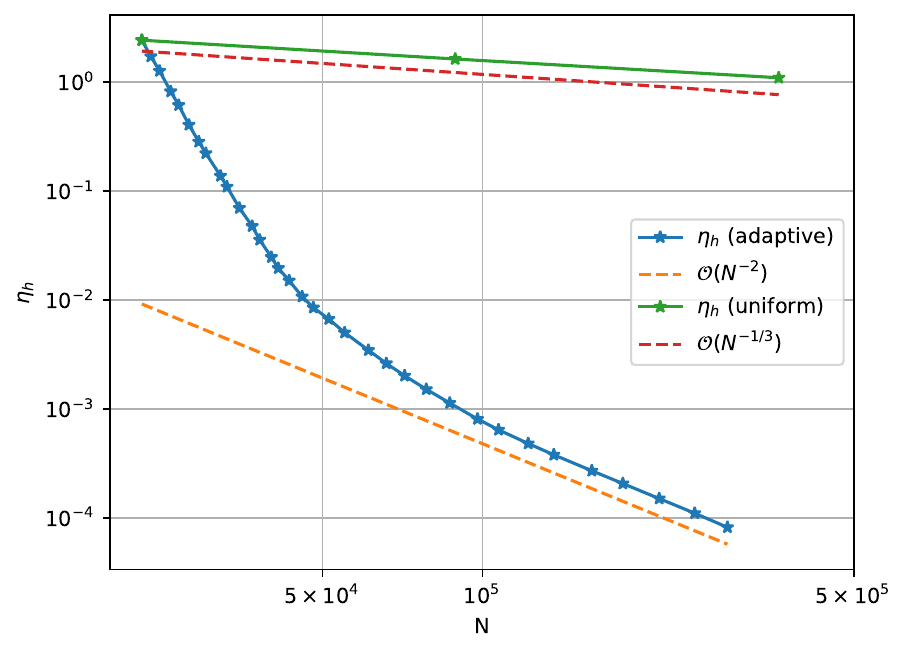}
\caption{Cook's membrane: convergence history of $\eta_h$ under adaptive and uniform refinement, AFW method with $k=4$, $\nu=0.3$ (left) and $\nu=0.49999$ (right).}
\label{fig:cook-conv}
\end{figure}

\begin{figure}[h!]
\centering
 \includegraphics[width=\textwidth]{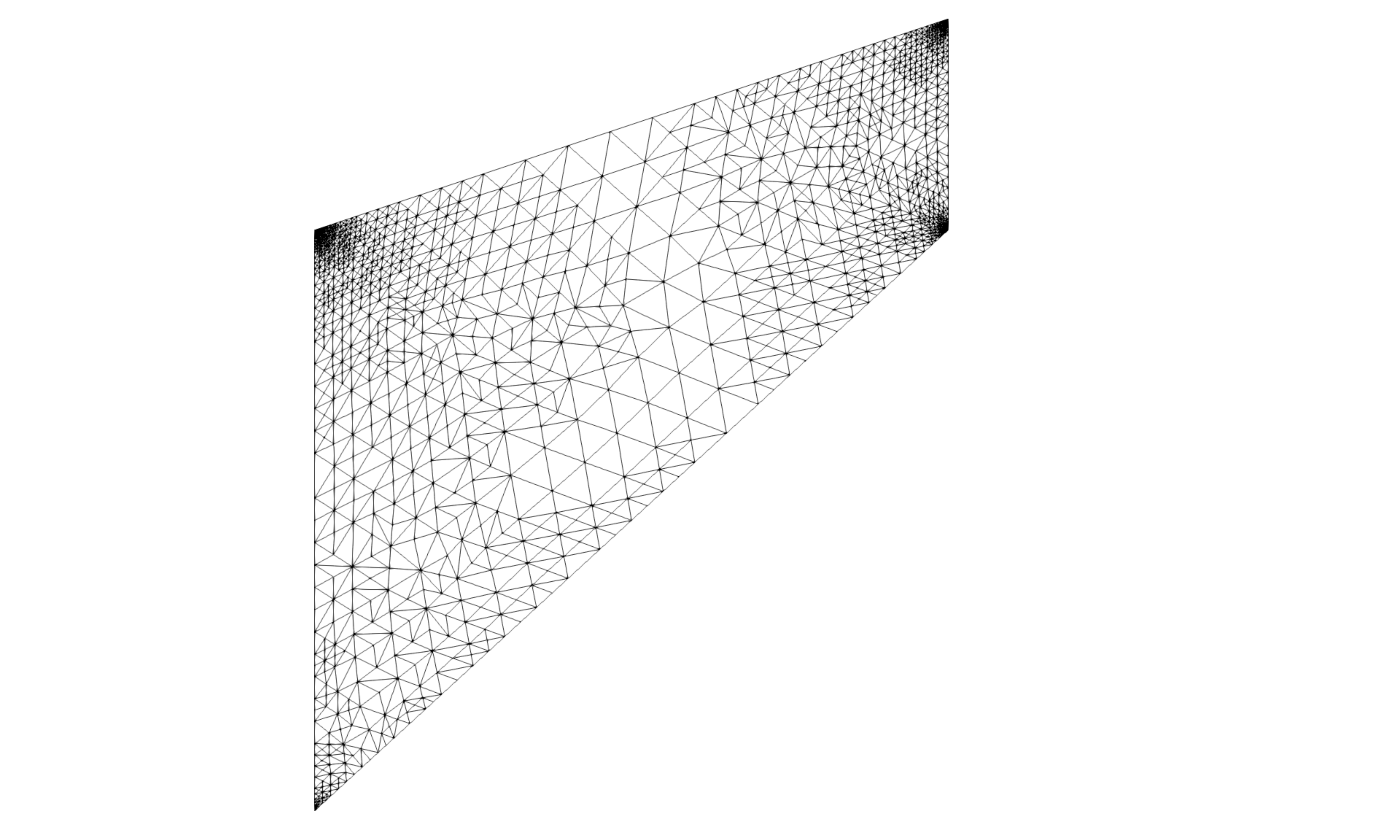}
\caption{Cook's membrane: adaptively refined mesh after 34 refinements, AFW method, $\nu=0.49999$.}
\label{fig:cook-mesh}
\end{figure}

Figure~\ref{fig:cook-conv} shows the convergence history of $\eta_h$ under adaptive versus uniform refinement: uniform refinement yields a markedly reduced rate, owing to the corner singularities, whereas adaptive refinement recovers the optimal rate $\CO(N^{-k/2})$. Figure~\ref{fig:cook-mesh} shows the resulting adaptively refined mesh, with refinement concentrated near the corners $(48,44)$ and $(48,60)$, where the traction data are incompatible, and near the corner $(0,44)$.

\begin{figure}[h!]
\centering
 \includegraphics[width=0.75\textwidth]{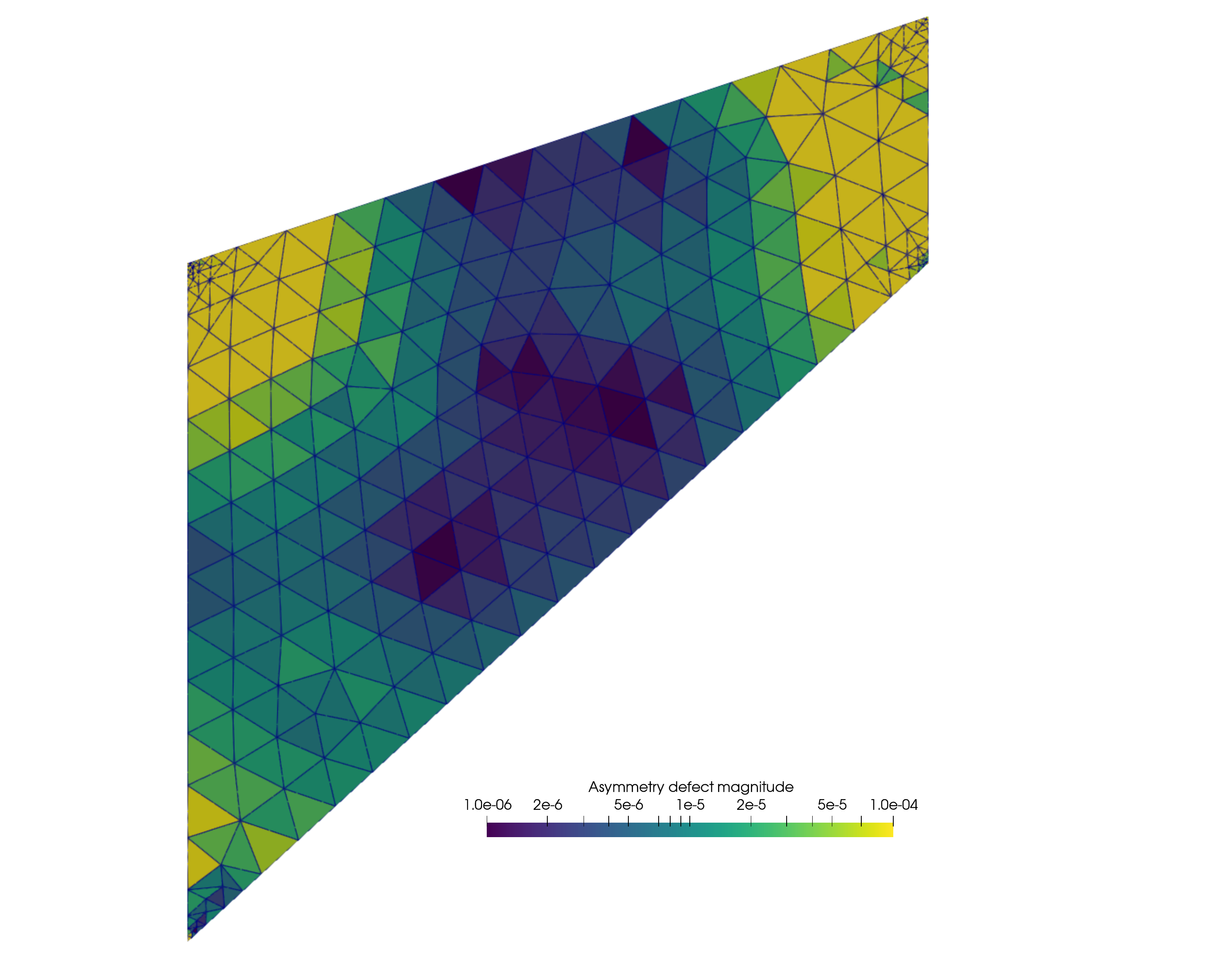}
\caption{Cook's membrane: elementwise asymmetry defect $\mu^{-1/2}\|\bsigma_h-\bsigma_h^{\T}\|_K$ after 15 adaptive refinements, AFW method, $k=4$, $\nu=0.49999$ (logarithmic color scale).}
\label{fig:cook-asym}
\end{figure}

As noted in Section~\ref{sec:aposteriori}, the local indicator $\eta_K$ carries the asymmetry-defect term $\mu^{-1}\|\bsigma_h-\bsigma_h^{\T}\|_K^2$, which measures, element by element, how far the discrete stress is from satisfying angular momentum balance. Figure~\ref{fig:cook-asym} shows its elementwise distribution, $\mu^{-1/2}\|\bsigma_h-\bsigma_h^{\T}\|_K$, after 15 adaptive refinements: the defect is largest around the corner $(0,44)$, where the clamped and traction-free parts of the boundary meet, and next to the loaded edge between $(48,44)$ and $(48,60)$, while it is smallest in the interior of the domain. This identifies where the asymmetry contributes most to the indicator.

\FloatBarrier

\section{Conclusions}
\label{sec:conclusions}

We have developed an a posteriori error estimator for weakly symmetric mixed finite element discretizations of linear elasticity, based on a local residual minimization problem obtained from a Stenberg-type postprocessing scheme, driven by the full displacement gradient rather than its symmetric part alone, by enlarging its test space. The resulting estimator is locally efficient, with the explicit constant $\sqrt{10}$, independent of $\lambda$; it is reliable in the compressible regime, with an explicit constant, and admits an alternative reliability estimate, with a constant that remains bounded as $\lambda\to\infty$.

As in the exactly symmetric case treated in a companion paper~\cite{CaceresVega}, a single indicator and a single comparison norm serve both regimes. The indicator carries, as do the estimators of~\cite{LS2024}, the asymmetry defect $\mu^{-1/2}\|\bsigma_h-\bsigma_h^{\T}\|$, quantifying the extent to which the discrete stress fails to satisfy angular momentum balance pointwise; its efficiency, as that of the remaining terms, follows from elementary arguments, without Korn-type inequalities or auxiliary bubble functions.

\section*{Acknowledgments}
The first author acknowledges the support provided by the Dirección de Investigación y Desarrollo (DICYT) of Universidad de Santiago de Chile under project DICYT Postdoctorado 042632VR\_Postdoc, and previously by the Agencia Nacional de Investigación y Desarrollo (ANID) through project FONDECYT 11251691. The second author gratefully acknowledges the support from ANID through project FONDECYT 11251691.

\bibliographystyle{elsarticle-num} 
\bibliography{references}

\end{document}